\documentclass[paper,twocolumn,twoside]{geophysics}

\usepackage{timet}
\usepackage{graphicx}
\usepackage{multirow}
\usepackage{multicol}
\usepackage{lipsum}
\usepackage{timet}
\usepackage{epsfig}
\usepackage{amsmath}
\usepackage{amsfonts}
\usepackage{amssymb}
\usepackage{color}
\numberwithin{equation}{section}
\usepackage[overload]{empheq}
\usepackage{algorithm}
\usepackage{algorithmic}
\usepackage{caption}
\usepackage{float}
\usepackage{graphics}
\usepackage{xcolor,graphicx}
\usepackage{csquotes}
\usepackage{mathtools}

\graphicspath{{Fig/}}

\begin{document}
\title{Robust Wavefield Inversion via Phase Retrieval}

\address{ \footnotemark[1]University of Tehran, Institute of Geophysics, Tehran, Iran, email: h.aghamiry@ut.ac.ir, agholami@ut.ac.ir \\ 
\footnotemark[2]Universit\'e C\^ote d'Azur, CNRS, Observatoire de la C\^ote d'Azur, IRD , G\'eoazur, Valbonne, France, email: aghamiry@geoazur.unice.fr, operto@geoazur.unice.fr}
\author{Hossein S. Aghamiry \footnotemark[1]\footnotemark[2], Ali Gholami \footnotemark[1] and St\'ephane Operto \footnotemark[2]}

\lefthead{A PREPRINT~~~~~~~~~~~~~~~~~~~~~~~~~~~~~~~~~~~~~~~~~~~~~~~~~~~~~~~~~~~~~~~~~~~~~~~~~~~~~~~~~~Aghamiry et al.}
\righthead{A PREPRINT~~~~~~~~~~~~~~~~~~~~~~~~~~~~~~~~~~~~~~~~~~~~~~~~~~~~~~~~~~~~~BTT regularized WIPR}

\maketitle

\begin{abstract}
\small{
Extended formulation of Full Waveform Inversion (FWI), called Wavefield Reconstruction Inversion (WRI), offers potential benefits of decreasing the nonlinearity of the inverse problem by replacing the explicit inverse of the ill-conditioned wave-equation operator of classical FWI (the oscillating Green functions) with a suitably defined data-driven regularized inverse. This regularization relaxes the wave-equation constraint to reconstruct wavefields that match the data, hence mitigating the risk of cycle skipping. The subsurface model parameters are then updated in a direction that reduces these constraint violations. However, in the case of a rough initial model, the phase errors in the reconstructed wavefields may trap the waveform inversion in a local minimum leading to inaccurate subsurface models. In this paper, in order to avoid matching such incorrect phase information during the early WRI iterations, we design a new cost function based upon phase retrieval, namely a process which seeks to reconstruct a signal from the amplitude of linear measurements. This new formulation, called Wavefield Inversion with Phase Retrieval (WIPR), further improves the robustness of the parameter estimation subproblem by a suitable phase correction.  We implement the resulting WIPR problem with an alternating-direction approach, which combines the Majorization-Minimization (MM) algorithm to linearise the phase-retrieval term and a variable splitting technique based upon the alternating direction method of multipliers (ADMM). This new workflow equipped with Tikhonov-total variation (TT) regularization, which is the combination of second-order Tikhonov and total variation regularizations and bound constraints, successfully reconstructs the 2004 BP salt model from a sparse fixed-spread acquisition using a 3~Hz starting frequency and a homogeneous initial velocity model. 
}
\end{abstract}


\section{INTRODUCTION}

Full waveform inversion (FWI) is a waveform matching procedure which provides a subsurface model with a resolution reaching to the half of the smallest propagated wavelength \citep{Virieux_2009_OFW}. 
Ultra-long offset wide-azimuth stationary-recording acquisitions provide a varied angular illumination of the subsurface, which is amenable to the development of the broadband velocity model \citep{Sirgue_2004_EWI}. 
However, the large number of propagated wavelengths generated by these acquisitions makes FWI prone to cycle skipping or phase wrapping \citep{Shah_PhD_2014, Choi_2015_UPI}. 
Various continuation strategies in frequencies, travel-times and offsets can be used to mitigate this pathology \citep[e.g., ][]{Shipp_2002_TDF,Gorszczyk_2017_TRW}. However, these multi-scale approaches can be tedious to implement and still require quite accurate initial velocity models and low frequencies to prevent cycle skipping at long offsets. 

To increase the resilience of FWI to cycle skipping, frequency-domain wavefield reconstruction inversion (WRI) has been proposed by \citet{VanLeeuwen_2013_MLM,vanLeeuwen_2016_PMP}. The governing idea of the method relies on the following statement: if we would be able to record a monochromatic wavefield generated by a single source everywhere in the subsurface, then it would be straightforward to reconstruct the subsurface model exactly in virtue of the wave equation bilinearity (Fig.~\ref{fig:mar_true}). Since we cannot record wavefields everywhere, the best we can do to approach the true wavefields is to find those that best match the sparse observations and satisfy the wave equation in a least-squares sense. This least-squares problem can be implemented by solving an overdetermined linear system gathering the weighted wave equation and the observation equation (the equation relating the data to the wavefield through a sparse sampling operator). The observation equation generates a wave equation relaxation whose strength depends on the accuracy of the available subsurface model and the weight assigned to the wave equation. Then, the subsurface parameters are updated from the reconstructed wavefields by the least-squares minimization of the wave equation errors (source residuals), and this two-step cycle is iterated until both the observation equation and the wave equation are satisfied with a prescribed accuracy. Typically, a small weight is assigned to the wave equation during early iterations to foster data fitting and prevent cycle skipping accordingly. This weight is progressively increased to guarantee that the wave equation is fulfilled near the convergence point. The parameter-estimation subproblem can be solved through a variable projection approach \citep{vanLeeuwen_2016_PMP} or through alternating optimization \citep{VanLeeuwen_2013_MLM}. In this latter case, the parameter-estimation subproblem is linearized around the reconstructed wavefield according to the wave-equation bilinearity.

WRI has been originally implemented with a penalty method, which requires a tedious dynamic control of the penalty parameter as above mentioned \citep{Fu_2017_DPM}. To overcome this issue, \citet{Aghamiry_2019_IWR} proposed to replace the penalty method with an augmented Lagrangian method equipped with operator splitting, namely the alternating-direction method of multiplier (ADMM) \citep{Boyd_2011_DOS}. Contrary to penalty methods, augmented Lagrangian methods converge to an accurate solution with a fixed penalty parameter by iteratively updating the Lagrange multipliers with the running sum of the constraint violations \citep[][Chapter 17]{Nocedal_2006_NO}.  This defines what is sometimes referred to as an iterative refinement or defect correction algorithm. Accordingly, \citet{Aghamiry_2019_IWR} called their method iteratively-refined (IR-) WRI.

ADMM also provides a suitable framework to implement bound constraints and  $\ell 1$-based nonsmooth regularizations (such as total variation (TV) regularization) in IR-WRI \citep{Aghamiry_2019_IBC}. IR-WRI with TV regularization was further improved by using more versatile Tikhonov-TV (TT) regularization combining second-order Tikhonov and TV regularizers through infimal convolution \citep{Aghamiry_2019_CRO}.
Infimal convolution means that the TT regularization explicitly decomposes the model into two components of different statistical properties (a smooth one and a blocky one) such that a suitable regularization can be tailored to each component (Tikhonov and TV) \citep{Gholami_2013_BCT}. 

Although IR-WRI, when equipped with basic frequency continuation strategies and advanced regularization techniques, mitigates cycle skipping and helps to overcome nonlinearity issues broadly speaking, still it can be trapped in a local minimum when the initial model is very far from the true model. 
In this case, the phase and amplitude of the reconstructed wavefields match those of the true wavefields only in the vicinity of the receivers, while they can be quite inaccurate elsewhere. In IR-WRI, the inaccuracies of the reconstructed wavefields are directly mapped into the subsurface model through a linear deconvolution-like processing. This implies that only the shallow part of the subsurface model is expected to be reliably updated during the early iterations of IR-WRI for surface acquisitions. It remains unclear from a theoretical viewpoint how these inaccuracies are progressively canceled out by IR-WRI in iterations through a depth-continuation process.

It is well acknowledged that the phase often has a more dominant role than amplitude in image processing and waveform inversion methods \citep{Oppenheim_1981_IPS, Shechtman_2015_PRW}. Therefore, the inaccuracies in the phase of the reconstructed wavefields may be those which can more likely drive the inversion toward spurious local minima. 
Accordingly, the objective of this study is to assess whether a reformulation of the parameter-estimation subproblem that mitigates the role of the phase during the early stages of the IR-WRI at the benefit of the more robust amplitude counterpart contributes to stabilizing the inversion.

To achieve this goal, we recast the parameter-estimation subproblem of IR-WRI as a phase retrieval problem \citep{Fienup_1982_PRA}. 
Phase retrieval considers the fundamental problem of how to reconstruct a signal from the magnitude of linear measurements (e.g., the magnitude of its Fourier transform). 
The name phase retrieval arose because if the phase of the linear measurements can be retrieved, then it is easy to reconstruct the signal by solving a linear problem \citep{Waldspurger_2015_PRM}.
It has a long history of applications in science and engineering when the phase of the linear measurements cannot be recorded or is inaccurate., e.g., optics \citep{Walther_1963_TQO}, X-ray crystallography \citep{Millane_1990_PRI, Harrison_1993_PPI}, astronomy \citep{Fienup_1982_PRA}. Also, it has recently attracted renewed interest in many imaging problems \citep{Fogel_2016_PRI}, such as  X-ray medical imaging \citep{Pfeiffer_2006_PRD,Burvall_2011_PRX}, optical imaging \citep{Shechtman_2015_PRW} and seismic processing \citep{Gholami_2014_PRT}.

In IR-WRI, the linear operator of the parameter estimation subproblem is formed by the so-called virtual sources \citep{Pratt_1998_GNF}, while the right-hand sides (the linear measurements) depend on the source and the Laplacian of the wavefields (considering the Helmholtz equation as wave equation).
According to the above definition of phase retrieval, we update the subsurface parameters from the magnitude of the right-hand sides rather than from their phase and amplitude. 
A potential downside of the phase-retrieval algorithm is ill-posedness generated by the rapidly-decreasing sensitivity of the inversion with depth, in particular in large contrast or attenuating media when a small amount of the seismic energy is transmitted in the deep subsurface. This makes the use of efficient regularization necessary, as we will show.
As in IR-WRI, we implement the iterative refinement procedure in the phase retrieval algorithm to design an adaptive control on the regularization parameters and achieve a faster convergence to the desired model. Once the velocity model estimation has been stabilized by phase retrieval during a first low-frequency band inversion, we switch to classical IR-WRI to process the higher frequencies with the phase information. 

We assess our phase-retrieval based waveform inversion, referred to as WIPR, against two complex synthetic case studies, the Marmousi II and the large contrast 2004 BP salt models with a sparse long-offset fixed-spread acquisition. We show how WIPR helps to reconstruct more accurately the phase of the reconstructed wavefields, which is translated into more accurate velocity estimation. Then, we show the impact of this more accurate early-stage velocity reconstruction on the convergence of IR-WRI at higher frequencies.

This paper is organized in a method, numerical, and discussion sections. In the method section, we first introduce phase retrieval in a general context. Second, we review the principles of WRI with iterative refinement and show how we implement phase retrieval in the underlying parameter-estimation subproblem. Then, we briefly review how compound regularizations can be implemented in WIPR. In the numerical experiment section, we show the application on the Marmousi and 2004 BP salt models. In the discussion section, we further discuss the resolution power of amplitudes in WIPR and show its impact on the phase reconstruction.

\section{NOTATION}
The mathematical symbols adopted in this paper are as follows.
We use italics for scalar quantities, boldface lower-case letters for vectors, and boldface capital letters for matrices and tensors.
We use the superscript $T$ to denote the adjoint of an operator and $*$ to show the complex conjugate. We use diag$(\bullet)$ to show a square diagonal matrix which includes vector $\bullet$ on its main diagonal.
The $i$th component of the column vector $\bold{x}$ is shown by $\it{x}_i$. 
For a complex number $\it{x}=\it{a}e^{j\it{b}}$ with $j=\sqrt{-1}$, $|\it{x}|=\it{a}$ denotes the magnitude of $\it{x}$ and $\angle \it{x}=\it{b}$ denotes its phase. For the $n$-length column vectors $\bold{x}$ and $\bold{y}$, the dot product is defined by $\langle \bold{x},\bold{y}\rangle=\bold{x}^T\bold{y}=\sum_{i=1}^n\it{x}_i^*\it{y}_i$ and
their Hadamard product, denoted by $\bold{x}\circ \bold{y}$, is another vector made up of their component-wise products, i.e. $(\bold{x}\circ \bold{y})_i=\it{x}_i\it{y}_i$.
The $\ell_1$- and $\ell_2$-norms of $\bold{x}$ are, respectively, defined by
$\|\bold{x}\|_1=\sum_{i=1}^n|\it{x}_i|$ and $\|\bold{x}\|_2=\sqrt{\langle \bold{x},\bold{x}\rangle}=\sqrt{\sum_{i=1}^n|\it{x}_i|^2}$. 

\section{METHOD}
In this section, we briefly review the phase retrieval problem \citep{Gerchberg_1972_APA} and a simple majorization-minimization (MM) \citep{Lange_2016_MOA} algorithm to solve this non-convex problem. Then we review the frequency-domain FWI \citep{Pratt_1998_GNF}, its WRI alternative \citep{vanLeeuwen_2016_PMP}, and its improved version called iteratively refined WRI (IR-WRI) \citep{Aghamiry_2019_IWR}. Finally, we introduce phase retrieval in the parameter-estimation subproblem of IR-WRI. In the following of this study, we refer to IR-WRI with phase retrieval as WIPR. 
WIPR is solved efficiently with the alternating direction method of multipliers (ADMM) \citep{Boyd_2011_DOS} and MM algorithms, which provide a suitable framework to implement bound constraints and TT regularization in the parameter estimation subproblem and tackle large data sets.

\subsection{Phase retrieval} \label{SECA}
Phase retrieval refers to the reconstruction of a complex signal from the amplitude of linear measurements.
Consider the complex-valued linear problem $\bold{Lx}=\bold{y}$ with $\bold{L}\in \mathbb{C}^{m\times n}$,  $\bold{x}\in \mathbb{C}^{n\times 1}$, and  $\bold{y}\in \mathbb{C}^{m\times 1}$
and assume that only the amplitude of the right hand side, $|\bold{y}|$, can be measured or is reliable.
The corresponding phase-retrieval problem can be written as \citep{Gholami_2014_PRT} 
 \begin{equation} \label{c0}
 \min_{\bold{x}} ~~ f(\bold{x})=\min_{\bold{x}}~~\frac12 \|\bold{|Lx|}-\bold{|y|}\|_2^2. 
 \end{equation}
Problem \eqref{c0} is non-convex where this non-convexity finds its root in removing the phase of the right hand side ($|\bold{Lx}|=|\bold{y}|$ is non-unique and there are many $\bold{x}$ that can fit the magnitude of $\bold{y}$).\\  
Several algorithms have been proposed to solve the problem
ranging from the Gerchberg-Saxton algorithm \citep{Gerchberg_1972_APA} as a first algorithm to 
more modern optimization techniques \citep{Candes_2013_PEA,Netrapalli_2013_PRU,Eldar_2016_RAP}.    
In this paper, we present a simple algorithm based on the MM technique (Appendix \ref{Appa}) to find a local minimum of $f(\bold{x})$ iteratively. Based on the proposed algorithm, the minimizer of the convex surrogate function 
\begin{align} \label{quad0}
\bold{x}^{k+1}=\underset{\bold{x}}{\arg\min}~~\frac12 \|\bold{Lx}-\bold{|y|}e^{j\angle \bold{Lx}^k}\|_2^2. 
\end{align} 
converges to the minimizer of $f(\bold{x})$ as iteration number, $k$,  tends to infinity. 
The quadratic problem \ref{quad0} is the same as the least-squares norm of the original system $\bold{Lx=y}$, when the phase of $\bold{y}$ is replaced with the phase of $\bold{Lx}^k$, which is extracted from the solution of $k$'th iteration.

%
%
\subsection{Full-waveform inversion versus wavefield inversion}
The reduced formulation of frequency-domain FWI can be written as \citep{Pratt_1998_GNF,Plessix_2006_RAS}
\begin{equation} \label{pratt_FWI}
\min_{\bold{m}} ~~ \|\bold{Pu}(\bold{m})-\bold{d}\|_2^2,
\end{equation} 
where $\bold{d}$ denotes the recorded seismic data, $\bold{P}$ is the linear observation operator that samples the wavefield $\bold{u(m)}$ at the receiver positions, and $\bold{m}$ denotes the subsurface parameters. The wavefield $\bold{u}(\bold{m})$ is the solution of the wave equation
\begin{equation}
\bold{u}(\bold{m})=\bold{A}^{-1}(\bold{m})\bold{b},
\label{eqwave}
\end{equation}
where $\bold{b}$ is the source term and $\bold{A(m)} \in \mathbb{C}^{n\times n}$ is the discretized wave-equation operator.
In this study, the wave equation is the Helmholtz equation, whose operator $\bold{A(m)}$ is given by  
\begin{equation}   \label{A}
\bold{A(m)} =  \bold{\Delta} + \omega^2 \bold{C} \text{diag}(\bold{m})\bold{B},
\end{equation}
where $\omega$ is the angular frequency, $\bold{m}$ contains the squared slownesses, $\bold{\Delta}$ is the discretized Laplace operator, $\bold{C}$ introduces boundary conditions (e.g., sponge-like absorbing boundary conditions such as perfectly-matched layers \citep{Berenger_1994_PML}) and $\bold{B}$ spreads the \enquote{mass} term $\omega^2\bold{C} \text{diag}(\bold{m})$ over all the  coefficients of the stencil to improve its accuracy following an anti-lumped mass strategy \citep{Chen_2013_OFD}.\\   

A main drawback of the reduced formulation in equation \eqref{pratt_FWI} is that the objective function depends on the model parameters $\bold{m}$ through the oscillating inverse operator $\bold{A}^{-1}(\bold{m})$.  
This makes the inverse problem highly nonlinear, and hence prone to convergence to inaccurate minimizer when the initial $\bold{m}$ is not accurate enough.

An alternative objective function as proposed by  \citet{vanLeeuwen_2016_PMP} is defined as
\begin{equation}  \label{WRI_penalty}
 \underset{\bold{m}}{\min} ~~ \|\bold{A(m})\bold{u}(\bold{m})-\bold{b}\|_2^2
\end{equation}
in which the wavefield is given by
\begin{equation} \label{u}
\bold{u}(\bold{m})=\left(\lambda\bold{A}(\bold{m})^T\bold{A}(\bold{m}) + \bold{P}^T\bold{P}\right)^{-1}
\left(\lambda\bold{A}(\bold{m})^T\bold{b} + \bold{P}^T\bold{d}\right).
\end{equation}
The wavefield $\bold{u}(\bold{m})$ is the solution of an augmented wave equation system, which gathers the wave equation $\bold{A(m})\bold{u}(\bold{m})=\bold{b}$ weighted by the scalar penalty parameter $\lambda>0$ and the observation equation $\bold{Pu}(\bold{m})=\bold{d}$ \citep[][ their equation 6]{VanLeeuwen_2013_MLM}.

A main advantage of equation \eqref{WRI_penalty} over equation \eqref{pratt_FWI} is that its objective function depends on the model parameters $\bold{m}$ through a regularized inverse of the wave-equation operator $\left(\lambda\bold{A}(\bold{m})^T\bold{A}(\bold{m}) + \bold{P}^T\bold{P}\right)^{-1}$, where the regularizer injects the prior information on the true wavefields, namely their sparse measurements $\bold{d}$, through the sampling operator $\bold{P}$. This wavefield reconstruction driven by the observations introduces a relaxation of the wave equation, namely source residuals, which are minimized to update the model parameters, equation~\ref{WRI_penalty}, such that the wavefields, equation~\ref{u}, are pushed back toward the wave equation.

The minimization problem, equation \eqref{WRI_penalty}, can be solved with a variable projection approach by enforcing the closed-form expression of $\bold{u}$, equation~\ref{u}, in the objective function, equation~\ref{WRI_penalty} \citep{vanLeeuwen_2016_PMP}, or with an alternating-direction strategy for $\bold{u}$ and $\bold{m}$ through a Gauss-Seidel like iteration \citep{VanLeeuwen_2013_MLM}. We will follow the latter option, where the model estimation problem, equation \eqref{WRI_penalty}, reduces to a quadratic optimization problem
\begin{equation}  \label{WRI}
\bold{m}^{k+1}= \underset{\bold{m}}{\arg\min} ~~ \|\bold{A(m})\bold{u}^k-\bold{b}\|_2^2,
\end{equation}
which, in virtue of the wave-equation bilinearity \citep{VanLeeuwen_2013_MLM,Aghamiry_2019_IWR}, can be recast as
\begin{equation}  \label{WRIb}
\bold{m}^{k+1}= \underset{\bold{m}}{\min} ~~ \|\bold{L}(\bold{u}^k)\bold{m}-\bold{y}(\bold{u}^k)\|_2^2,
\end{equation}
with
\begin{equation} \label{L}
\begin{cases}
 \bold{L}(\bold{u}^k) = \omega^2  \bold{C}\text{diag}(\bold{B}\bold{u}^k),\\
 \bold{y}(\bold{u}^k)=\bold{b}-\Delta \bold{u}^k.
 \end{cases} 
\end{equation}
and $\bold{u}^k \equiv \bold{u}(\bold{m}^k)$. 

During wavefield reconstruction, equation~\ref{u}, a small value of $\lambda$ is used during the early iterations to foster the data fit, hence mitigating cycle skipping, at the expense of the fidelity with which the wave equation is satisfied. Then, $\lambda$  is progressively increased in iterations such that the wave equation is satisfied
at the convergence point \citep{Fu_2017_DPM}.
To avoid the tedious dynamic control of penalty parameters and converge in an automatic way toward accurate solution with fixed $\lambda$, \citet{Aghamiry_2019_IWR} introduced iteratively-refined WRI (IR-WRI). IR-WRI relies on an iterative refinement procedure, which finds its root in augmented Lagrangian method \citep[][ Chapter 17]{Nocedal_2006_NO}. In the end, it simply consists of updating at the end of each iteration the right-hand sides of the wave equation and the observation equation, namely the data $\bold{d}$ and the source $\bold{b}$, with the running sum of the data and source residuals in iteration
\begin{subequations}
\begin{align}
\bold{b}^{k+1} &= \bold{b}^k +\bold{b} - \bold{A(m}^{k+1})\bold{u}^{k+1}, \label{IR_WIPR_b}\\
\bold{d}^{k+1} &= \bold{d}^k +\bold{d} - \bold{Pu}^{k+1}. \label{IR_WIPR_d} 
\end{align}
\end{subequations}
In the framework of linear inverse problems, $\bold{d}^{k+1}$ and  $\bold{b}^{k+1}$ record the  wavefield and model refinements performed at previous iterations to refine the wavefields and the model parameters at the current iteration from the residual errors only. This iterative refinement or defect correction \citep{Bohmer_1984_DCM} procedure, which is lacking in WRI, is the key ingredient which controls the accuracy of the wavefield and model solutions at the convergence point when fixed $\lambda$ is used. The reader is also referred to \citet {Gholami_2017_CNA, Gholami_2019_3DD} for other recent applications of iterative refinement in the field of geophysics.

IR-WRI has shown promising results even for complicated velocity models when equipped with frequency continuation strategies and appropriate regularization \citep{Aghamiry_2019_IBC, Aghamiry_2019_CRO}. 
However, when the initial model is far from the true model, the phase of the reconstructed wavefields becomes inaccurate in areas that are located far away from the receivers. These phase inaccuracies can lead to an inaccurate model reconstruction of these regions, hence trapping the inversion into a local minimum.
In order to mitigate this issue, we show in the next section how to solve the quadratic optimization problem, equation~\ref{WRIb}, with phase retrieval.

\subsection{Wavefield inversion with phase retrieval (WIPR)}
To estimate the model parameters without involving inaccurate phase information, we propose to replace the IR-WRI parameter-estimation problem given by equation \ref{WRIb} by the following phase retrieval problem 
\begin{equation}  \label{WIPR_penalty}
\bold{m}^{k+1} = \underset{\bold{m}}{\arg\min} ~~ \||\bold{L}(\bold{u}^k)\bold{m}|-|\bold{y}(\bold{u}^k)|\|_2^2.
\end{equation}  
In order to solve \eqref{WIPR_penalty}, we replace it by a surrogate majorizing function (as described in section \ref{SECA} and Appendix \ref{Appa})
\begin{equation}  \label{mains}
\bold{m}^{k+1} = \underset{\bold{m}}{\arg\min} ~~ 
 \|\bold{L}(\bold{u}^k)\bold{m}-|\bold{y}(\bold{u}^k)| e^{j\angle \bold{L}(\bold{u}^k)\bold{m}^k}\|_2^2.
\end{equation}
This objective function is now quadratic and admits a closed form solution.

In order to better understand why phase retrieval helps to improve the solution, we give a tentative interpretation hereafter. From equation \eqref{u} it is seen that, at each iteration, the wavefield approximately satisfies the wave equation
\begin{equation}
\bold{A(m}^k)\bold{u}^{k}\approx \bold{b}
\end{equation}
or
\begin{equation} \label{33}
\bold{L}(\bold{u}^{k})\bold{m}^k\approx \bold{y}(\bold{u}^{k}),
\end{equation}
where the approximation level is controlled by the parameter $\lambda$ and the approximation  becomes an equality as $\lambda\rightarrow \infty$.
The next iterate $\bold{m}^{k+1}$ is then found such that it minimizes the wave-equation errors generated by the wave-equation relaxation, equation \eqref{WRI}. In this way, simple minimization via equation \eqref{WRI} forces the model to match both the amplitude and phase information of $\bold{y}(\bold{u}^{k})$. 
In order to reduce the imprint of the phase errors, a phase correction or phase alignment step can be applied to equation \eqref{33} before updating the model parameters:
\begin{equation} \label{ph}
\underset{\bold{\phi}}{\min}~ \|\bold{L}(\bold{u}^{k})\bold{m}^k-|\bold{y}(\bold{u}^{k})|e^{j\phi} \|_2^2.
\end{equation}
The optimal minimizer of equation \eqref{ph} is given by $\bold{\phi}=\angle\bold{L}(\bold{u}^{k})\bold{m}^k$.
Using the optimal $\bold{\phi}$, the next iterate is found by minimizing the corrected quadratic problem, as presented in equation \eqref{mains}.
In this way, the solution of the inverse problem is less affected by the phase errors in $\bold{u}$ and thus hopefully less prone to convergence to a local minimum.
The reader is referred to  \citet{Jiang_2017_RPR,Qian_2017_IAO}  for other applications of phase retrieval with alternating optimization. 

Furthermore, we apply the iterative refinement procedure to WIPR in the same way as for IR-WRI by updating $\bold{b}$ and $\bold{d}$ according to equations \ref{IR_WIPR_b} and \ref{IR_WIPR_d} as outlined in the Algorithm \ref{Alg_IR_WIPR}. 
This implies that, although we recast the parameter-estimation subproblem as a phase retrieval problem (Algorithm \ref{Alg_IR_WIPR}, Line 5), we use both the phase and amplitude of the source and data residuals to update the right-hand sides of the wavefield-reconstruction and parameter-estimation subproblems (Algorithm \ref{Alg_IR_WIPR}, Lines 8 and 9).
\begin{algorithm}[!h]
\caption{WIPR algorithm}  \label{Alg_IR_WIPR}
\begin{algorithmic}[1]
 \STATE Input: $\bold{m}^0$
\STATE initialize: $\bold{A}\leftarrow\bold{\Delta} + \omega^2 \bold{C} \text{diag}(\bold{m}^{0})\bold{B}$,
 ${\bold{b}^0}\leftarrow\bold{0}$, ${\bold{d}^0}\leftarrow\bold{0}$    \\
\WHILE{convergence conditions not satisfied}
\STATE $\bold{u}^{k+1}\leftarrow\left[\bold{P}^T\bold{P}+\lambda \bold{A}^T\bold{A}\right]^{-1}[\bold{P}^T(\bold{d}+{\bold{d}^k})+\lambda\bold{A}^T(\bold{b}+{\bold{b}^k}) ]$\\
\STATE $\bold{L} \leftarrow \omega^2  \bold{C}\text{diag}(\bold{B}\bold{u}^{k+1})$\\
$\tilde{\bold{y}}\leftarrow |\bold{b} +{\bold{b}^k}-\Delta\bold{u}^{k+1}|e^{j\angle \bold{L}\bold{m}^k}$\\
\STATE $\bold{m}^{k+1}\leftarrow \left[\bold{L}^T\bold{L}\right]^{-1}\bold{L}^T\tilde{\bold{y}}$\\
\STATE $\bold{A}\leftarrow\bold{\Delta} + \omega^2 \bold{C} \text{diag}(\bold{m}^{k+1})\bold{B}$\\
\STATE ${\bold{b}^{k+1}} \leftarrow {\bold{b}^{k}} +\bold{b} - \bold{A}\bold{u}^{k+1}$ \\
\STATE $\bold{d}^{k+1}\leftarrow \bold{d}^k +\bold{d} - \bold{Pu}^{k+1}$
\ENDWHILE
\end{algorithmic}
\end{algorithm}

\subsection{Tikhonov-total variation (TT) regularized WIPR with bounding constraints}
In this section, we implement TT regularization and bounding constraints in WIPR to stabilize the updates of the model parameters.
By adding a convex regularization term $\|\bold{m} \|_{\text{TT}}$ to the cost function of WIPR, we get the following optimization for the $\bold{m}$-subproblem
\begin{equation}  \label{WIPR_TT}
\bold{m}^{k+1}= \underset{\bold{m}\in \mathcal{C}}{\arg\min} ~~ \|\bold{m} \|_{\text{TT}} + \lambda
 \|\bold{L}(\bold{u}^k)\bold{m}-|\bold{y}(\bold{u}^k)| e^{j\angle \bold{L}(\bold{u}^k)\bold{m}^k}\|_2^2,
\end{equation}
where $\mathcal{C} = \{\bold{x} \in \mathbb{R}^{n\times 1}~\vert~ \bold{m}_{{l}} \leq \bold{x} \leq \bold{m}_{{u}}\}$ is the set of all feasible models bounded by the lower bound $\bold{m}_{{l}}$ and the upper bound $\bold{m}_{{u}}$ and $\|\bold{m} \|_{\text{TT}}$ is the TT regularization functional as defined by \citet{Aghamiry_2019_IBC,Aghamiry_2019_CRO}
\begin{equation}
 \label{TT_WIPR}
\|\bold{m} \|_{\text{TT}}=\underset{\bold{m}=\bold{m}_1+\bold{m}_2}{\min} \|\bold{m}_1 \|_{\text{TV}}+\alpha \|\bold{m}_2 \|_{\text{Tikh}},
\end{equation} 
with 
\begin{equation}
\begin{cases}
\|\bold{m}\|_{\text{TV}} =\sum  \sqrt{|\nabla_{\!\! x}\bold{m}|^2 + |\nabla_{\!\! z}\bold{m}|^2}, \\
\|\bold{m}\|_{\text{Tikh}} =\sum  \left(|\nabla_{\!\! xx}\bold{m}|^2 +2|\nabla_{\!\! xz}\bold{m}|^2+ |\nabla_{\!\! zz}\bold{m}|^2\right),
\end{cases}
\end{equation}
in which 
the sum runs over all element, $\nabla_{\!\! i}$ are first order difference operators in direction $i$ and
  $\nabla_{\!\! ij}$ are second-order differential operators in directions  $i$ and $j$. 

The TT regularization is an infimal convolution-based combination of the second-order Tikhonov and TV regularizations and is suitable for recovering piecewise-smooth models \citep{Gholami_2013_BCT, Aghamiry_2019_CRO}. It decomposes the model $\bold{m}$ explicitly into two components, $\bold{m}_1$ and $\bold{m}_2$, which have different statistical properties. Here, the blocky component $\bold{m}_1$ and the smooth component $\bold{m}_2$ are captured by the TV and the Tikhonov regularizations, respectively.

Subproblem \eqref{WIPR_TT} is solved with ADMM and the split-Bregman scheme, which de-couples the  $\ell_1$ and $\ell_2$ components and bound constraints of the function through the introduction of auxiliary variables and solves each related subproblem in sequence \citep{Goldstein_2009_SBM}. 
By doing this, we come up with a least-squares problem with a closed-form expression to update $\bold{m}_1$ and $\bold{m}_2$, a proximity subproblem to update the auxiliary variables and a gradient ascent step to update dual variables. The reader is referred to \citet{Aghamiry_2019_CRO} and \citet[Section 2.2.2]{Aghamiry_2019_IBC} for the detailed derivation of the algorithm.

\section{NUMERICAL EXAMPLES}
We assess our workflow against the Marmousi II model \citep{Martin_2006_M2E} and the large-contrast 2004 BP salt model \citep{Billette_2004_BPB} when the inversion is started from a homogeneous velocity model and a 3~Hz frequency such that a considerable phase error is generated.
In both cases, we consider a long-offset fixed-spread acquisition and perform the frequency-domain inversion in the 3~Hz - 13~Hz frequency band with a frequency interval of 0.5~Hz. Small batches of two frequencies with one frequency overlapping between two consecutive batches are successively inverted following a classical frequency continuation strategy. We apply the phase retrieval algorithm only during the first frequency batch processing to stabilize the inversion. Then, we proceed with regular IR-WRI at higher frequencies. The IR-WRI results obtained from an early WIPR stage will be referred to as IR-WRI$_{pr}$ to prevent confusions with results obtained when IR-WRI is applied from the first frequency-batch inversion. During IR-WRI/IR-WRI$_{pr}$, we perform three paths through the frequency batches to improve the inversion results, using the final model of one path as the initial model of the next one (these cycles can be viewed as outer iterations of the algorithms). The starting and finishing frequencies of the paths are [3.5, 6], [4, 8.5], [6, 13]~Hz, respectively, where the first element of each pair shows the starting frequency and the second one is the finishing frequency. We will compare the results of IR-WRI$_{pr}$ and IR-WRI to assess the improvements resulting from the initial WIPR step.
We will also compare the results of IR-WRI and IR-WRI$_{pr}$ when regularization is applied or not such that the relative role of phase retrieval and regularization can be discriminated. Namely, we will illustrate that the role of phase retrieval is to extend the linear regime of the waveform inversion further when an inaccurate starting model is used, while regularization injects suitable prior to manage non-uniqueness of the solution in the poorly illuminated area. When regularization is applied, the same tuning is used during WIPR and IR-WRI, and during IR-WRI and IR-WRI$_{pr}$. For all the numerical tests, we use a 9-point finite-difference staggered-grid stencil with PML boundary condition along the four edges of the model (no free-surface boundary condition is used) and anti-lumped mass to solve the Helmholtz equation, where the stencil coefficients are optimized to the frequency \citep{Chen_2013_OFD}. 

\subsection{Marmousi II model}
\begin{figure}[ht!]
\centering
\includegraphics[width=0.48\textwidth]{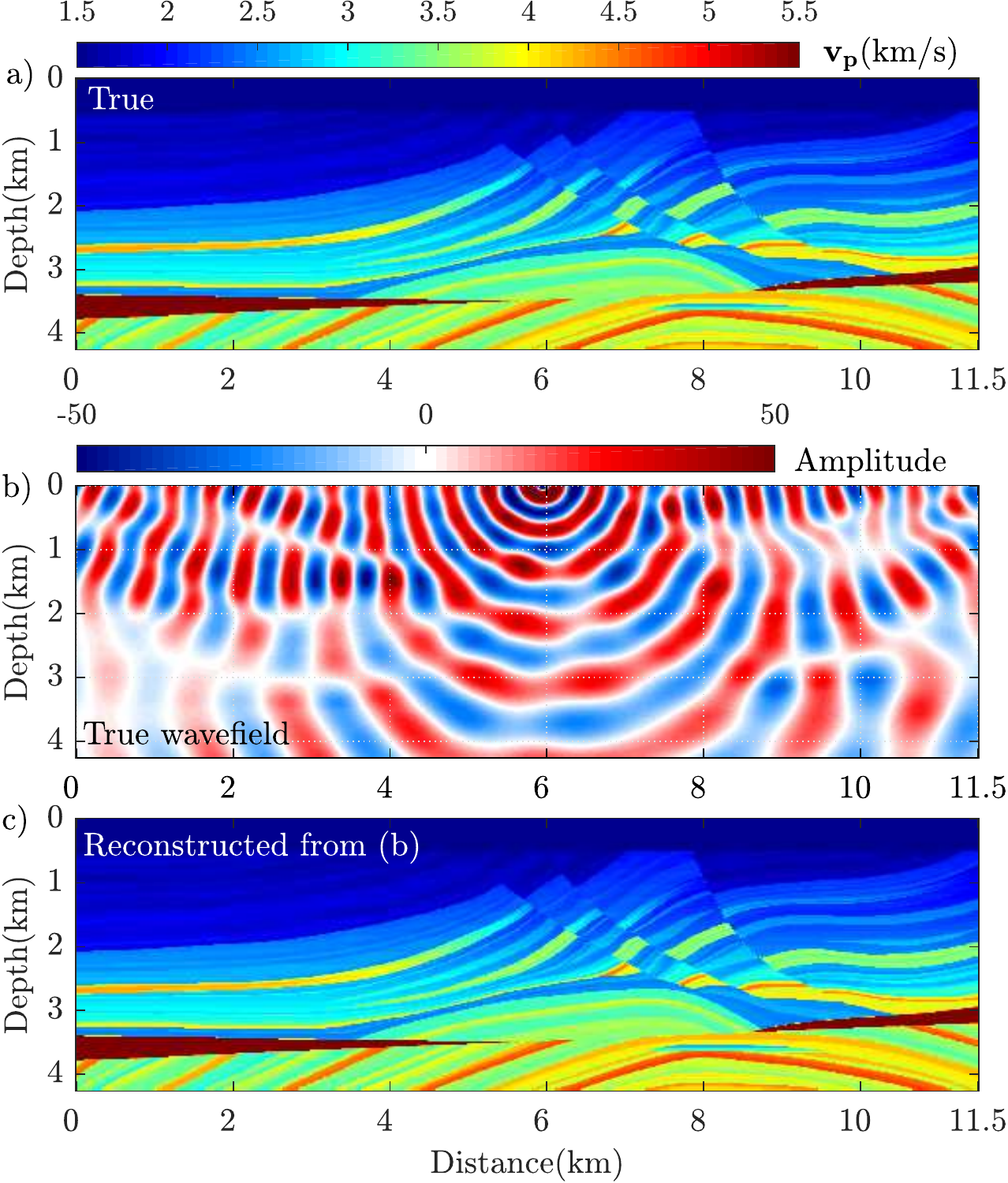}
\caption{Illustration of the wave-equation bilinearity. (a) True Marmousi II velocity model. (b) Three-Hz wavefield (real part) for the source located at the surface and a distance of 6 km. (c) Velocity model inferred from the amplitude and phase of the monochromatic wavefield, equation \ref{eqbilinpa}. The same model can be inferred from the sole amplitude of the monochromatic wavefield, equation \ref{eqbilina}. The reconstructed velocity model exactly matches the true velocity model. }
\label{fig:mar_true}
\end{figure}

The Marmousi II model covers a $11500~m \times 4200~m$ spatial domain (Fig.~\ref{fig:mar_true}a). The fixed-spread surface acquisition consists of 56 sources spaced 200~m apart, and 230 receivers spaced 50~m apart at the surface. The source signature is a 10~Hz Ricker wavelet, and the wavespeed in the homogeneous starting model is 3~km/s.

\begin{figure}[ht!]
\centering
\includegraphics[width=0.48\textwidth]{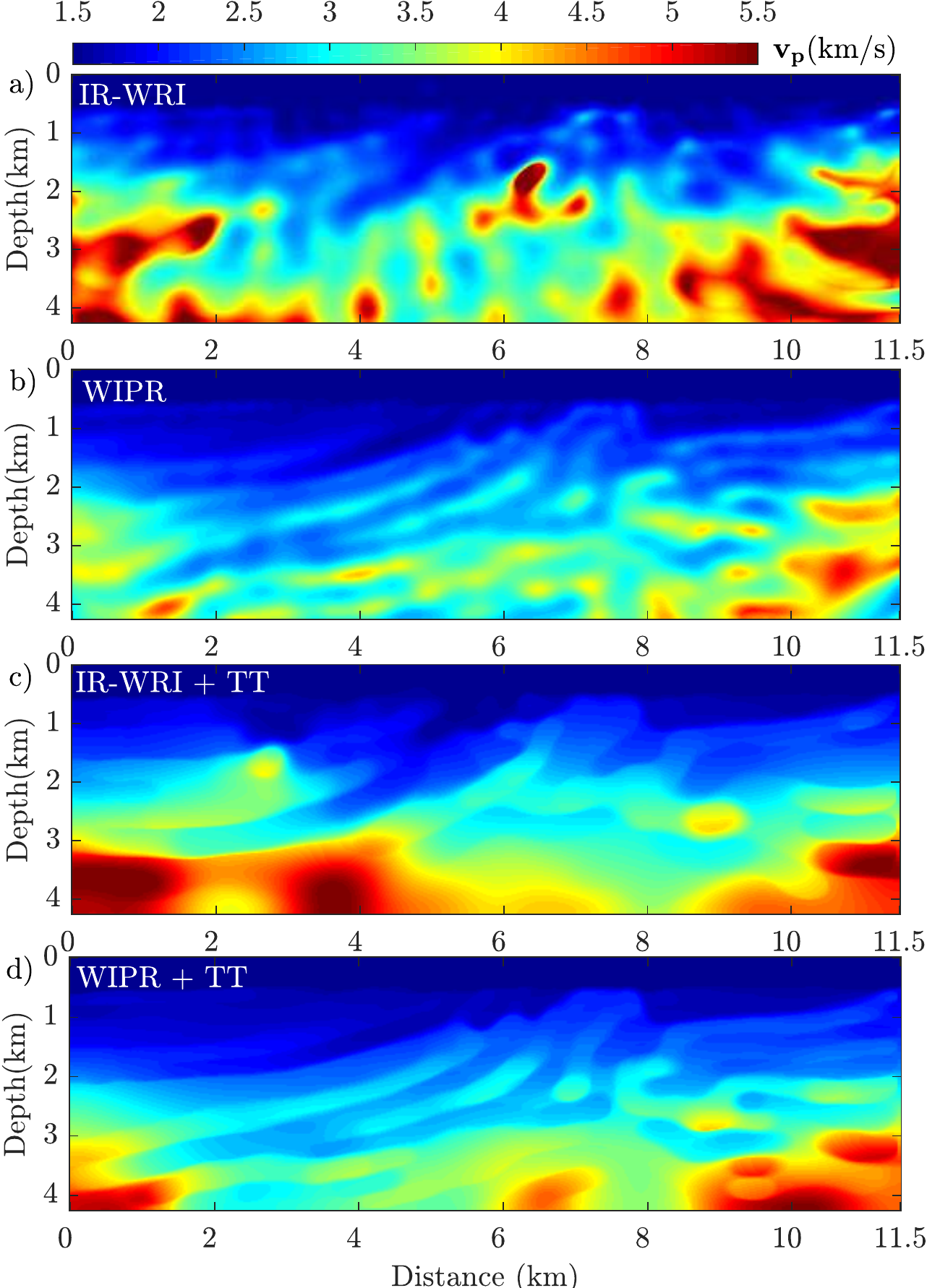}
\caption{Marmousi II example. Velocity models after the \{3,3.5\}~Hz inversion. (a) IR-WRI, (b) WIPR, (c) TT regularized IR-WRI, and (d) TT regularized WIPR. Bound constraints are used for all of the tests from the first iteration. }
\label{fig:mar_model_first}
\end{figure}

\begin{figure}[ht!]
\centering
\includegraphics[width=0.48\textwidth]{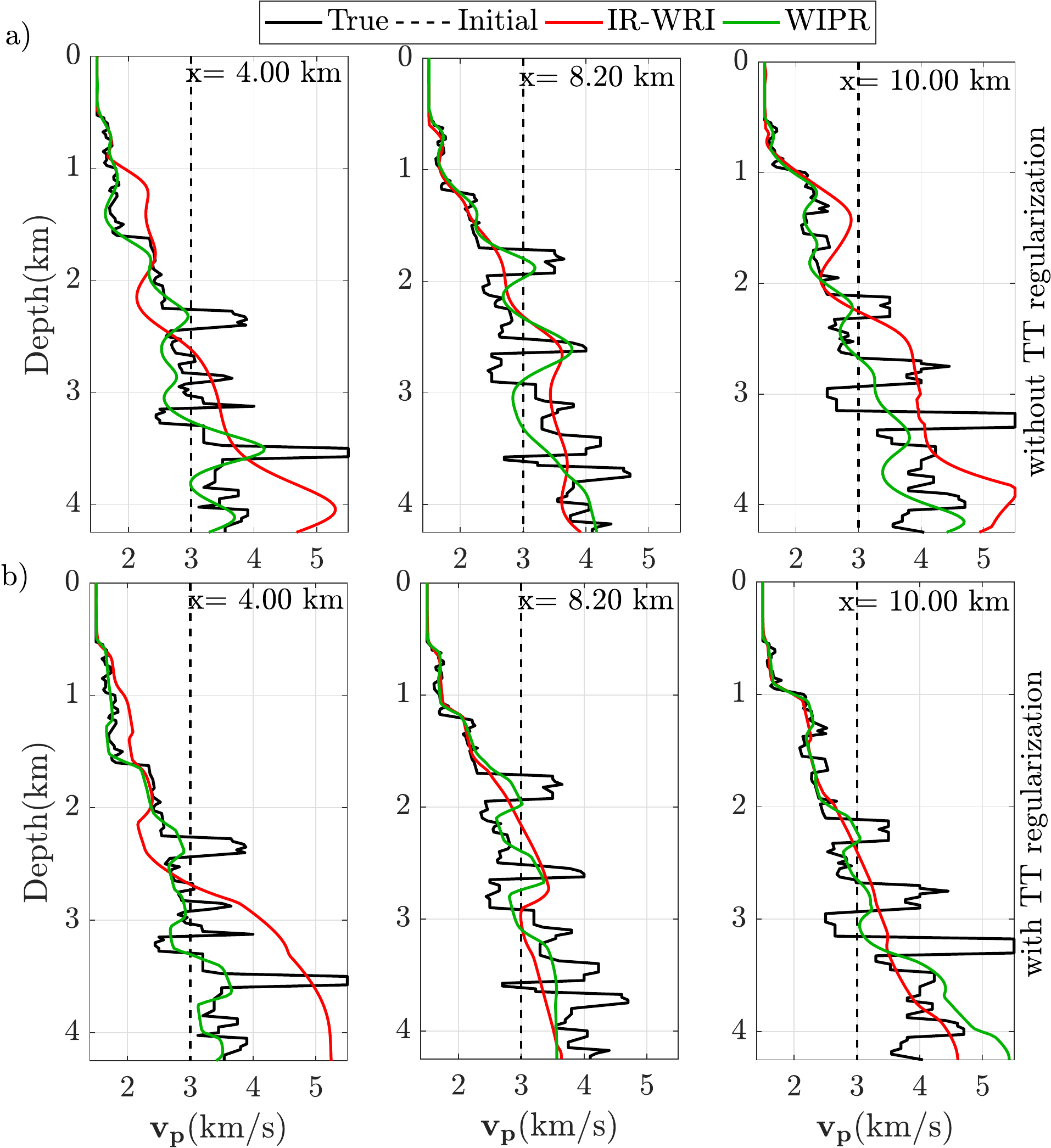}
\caption{Marmousi II example. Direct comparison along with three vertical logs at $x=4, 8.2$, and $10~km$ between the true velocity model (black), the initial model (dashed black), and the IR-WRI (red) and WIPR (green) models shown in Fig.~\ref{fig:mar_model_first}. (a) Without TT regularization. (b) With TT regularization.}
\label{fig:mar_log_first}
\end{figure}

We first compare the results of WIPR and IR-WRI after 45 iterations of the [3,3.5]~Hz frequency batch inversion when bound constraints are applied without any additional regularization (Fig. \ref{fig:mar_model_first}a-b). Bound constraints are used for all of the tests of this section from the first iteration. WIPR manages to reconstruct an accurate long-wavelength approximation of the Marmousi model (Fig.~\ref{fig:mar_model_first}b), while IR-WRI clear fails to capture the kinematic trend of the model (Fig.~\ref{fig:mar_model_first}a).
This statement is clearly illustrated by the direct comparison between the true model, the initial model, and the reconstructed IR-WRI and WIPR models along with three vertical logs at horizontal distances of 4, 8.2, and 10~km (Fig.~\ref{fig:mar_log_first}a). We also apply TT regularization during WIPR and IR-WRI (Figures \ref{fig:mar_model_first}c-d and \ref{fig:mar_log_first}b). The TT regularization improves the WIPR and IR-WRI results in the poorly-illuminated area, although it does not remove high-velocity artefacts near the left-bottom end of the IR-WRI velocity model (Fig.~\ref{fig:mar_log_first}b, x=4~km). \\

%
%
%
%
\begin{figure}[ht!]
\centering
\includegraphics[width=0.48\textwidth]{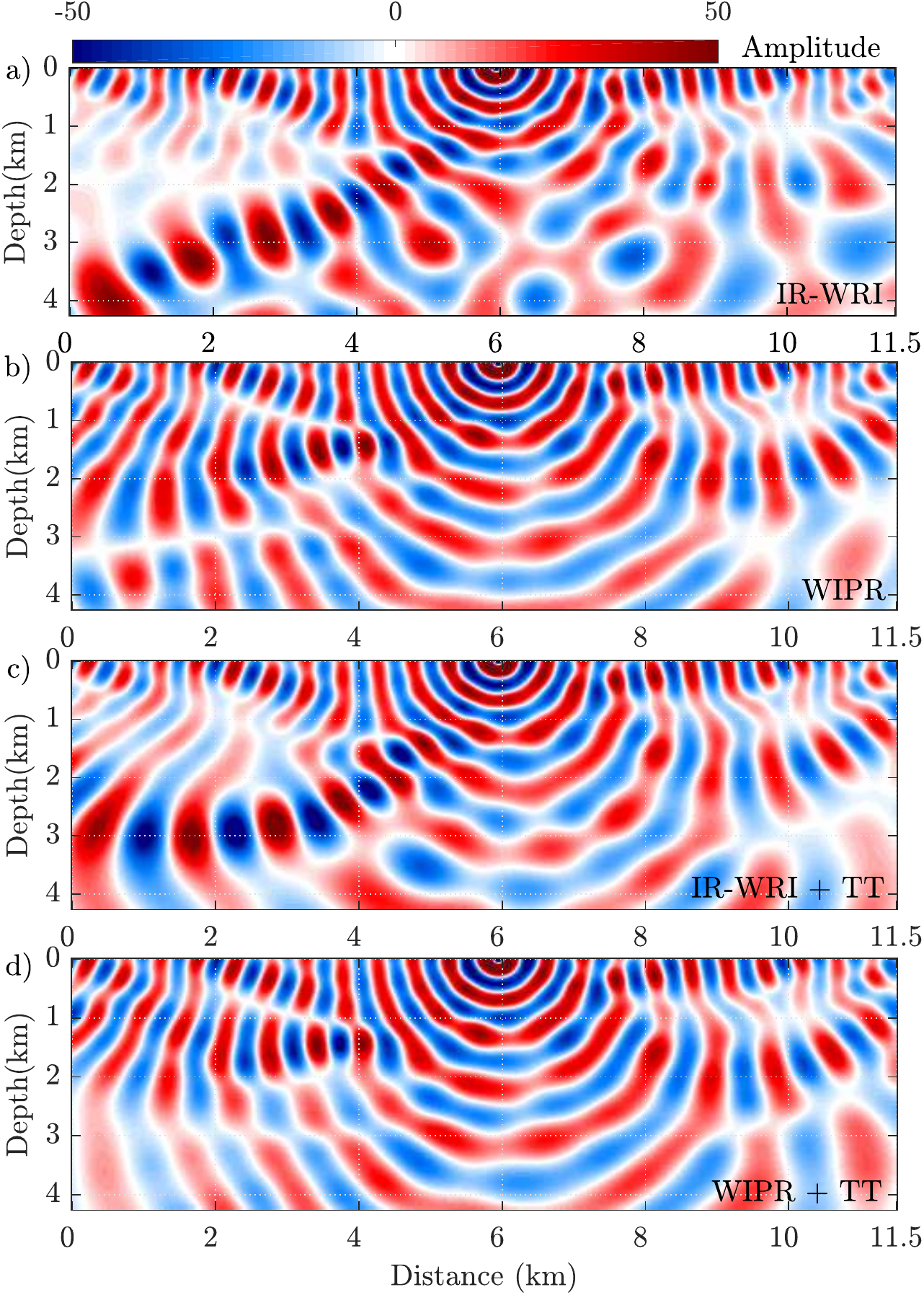}
\caption{Marmousi II example. Three-Hz wavefields (real part) reconstructed after the first frequency batch inversion for the source located at 6~km. (a) IR-WRI, (b) WIPR, (c) TT regularized IR-WRI, and (d) TT regularized WIPR. The corresponding velocity models are shown in Fig.~\ref{fig:mar_model_first}. These wavefields can be compared with the true one shown in Fig.~\ref{fig:mar_true}b.}
\label{fig:mar_wavefield_first}
\end{figure}

Since the model parameters are updated from the reconstructed wavefields, the accuracy of these wavefields is a visible indicator of the reliability of the estimated model parameters. We compare the real part of the 3~Hz wavefields reconstructed by IR-WRI and WIPR without and with TT regularization (Fig. \ref{fig:mar_wavefield_first}) with the true wavefield (Fig.~\ref{fig:mar_true}b) for a source located at the horizontal distance of 6~km. Both IR-WRI and WIPR wavefields are in phase with the true wavefield along the receiver line since we tune the penalty parameter such that the data are matched from the first iteration. However, the WIPR wavefield matches much better the true wavefield than the IR-WRI counterpart in depth. When regularization is used, the accuracy of the wavefield reconstruction is mostly improved in the poorly-illuminated zones, mainly near the bottom ends of the model. \\

\begin{figure}[ht!]
\centering
\includegraphics[width=0.48\textwidth]{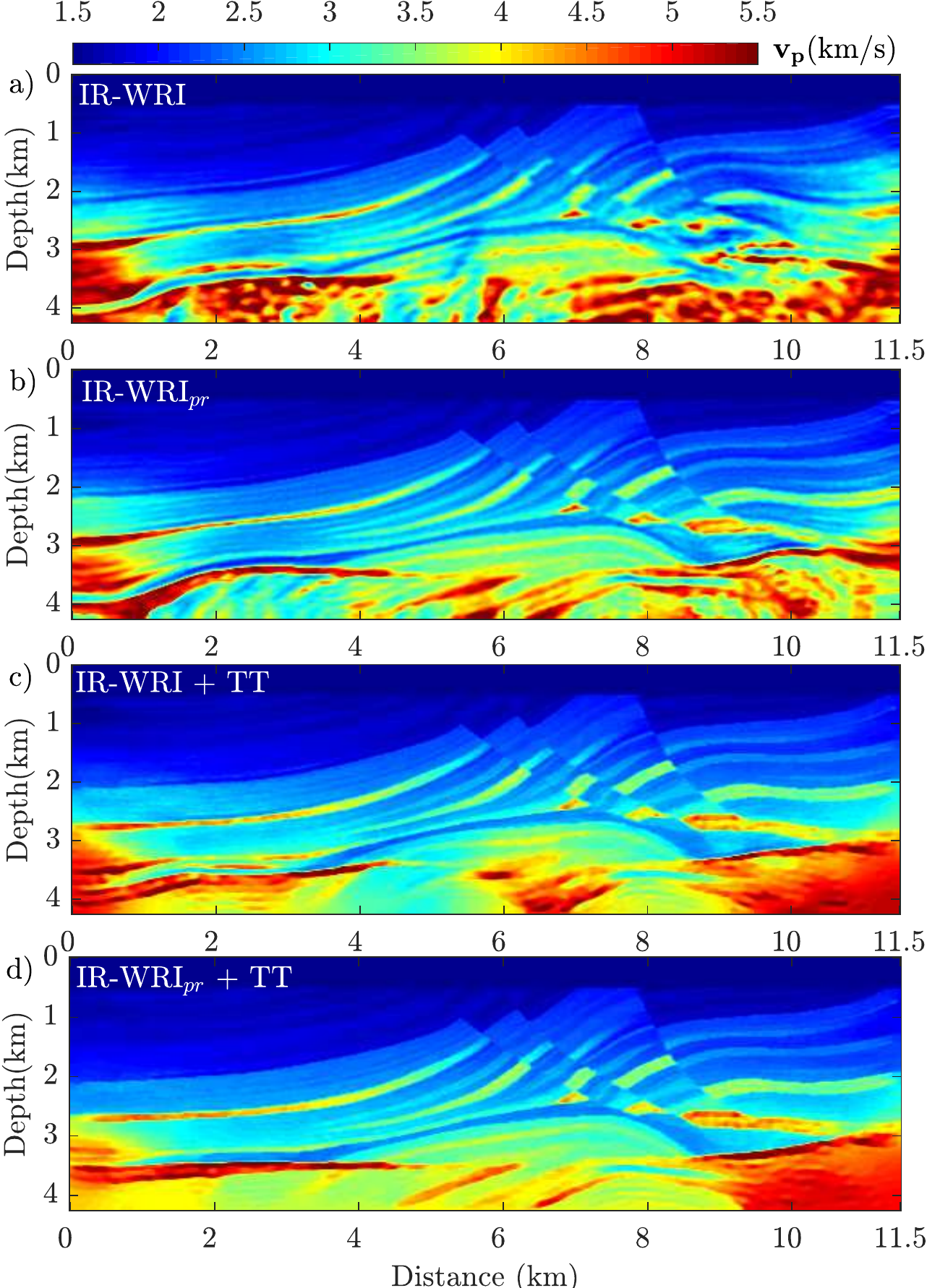}
\caption{Marmousi II example. Final velocity models obtained by IR-WRI (a,c) and IR-WRI$_{pr}$ (b,d). (a,b) Without TT regularization. (c,d) With TT regularization.}
\label{fig:Mar_final}
\end{figure}

\begin{figure}[ht!]
\centering
\includegraphics[width=0.48\textwidth]{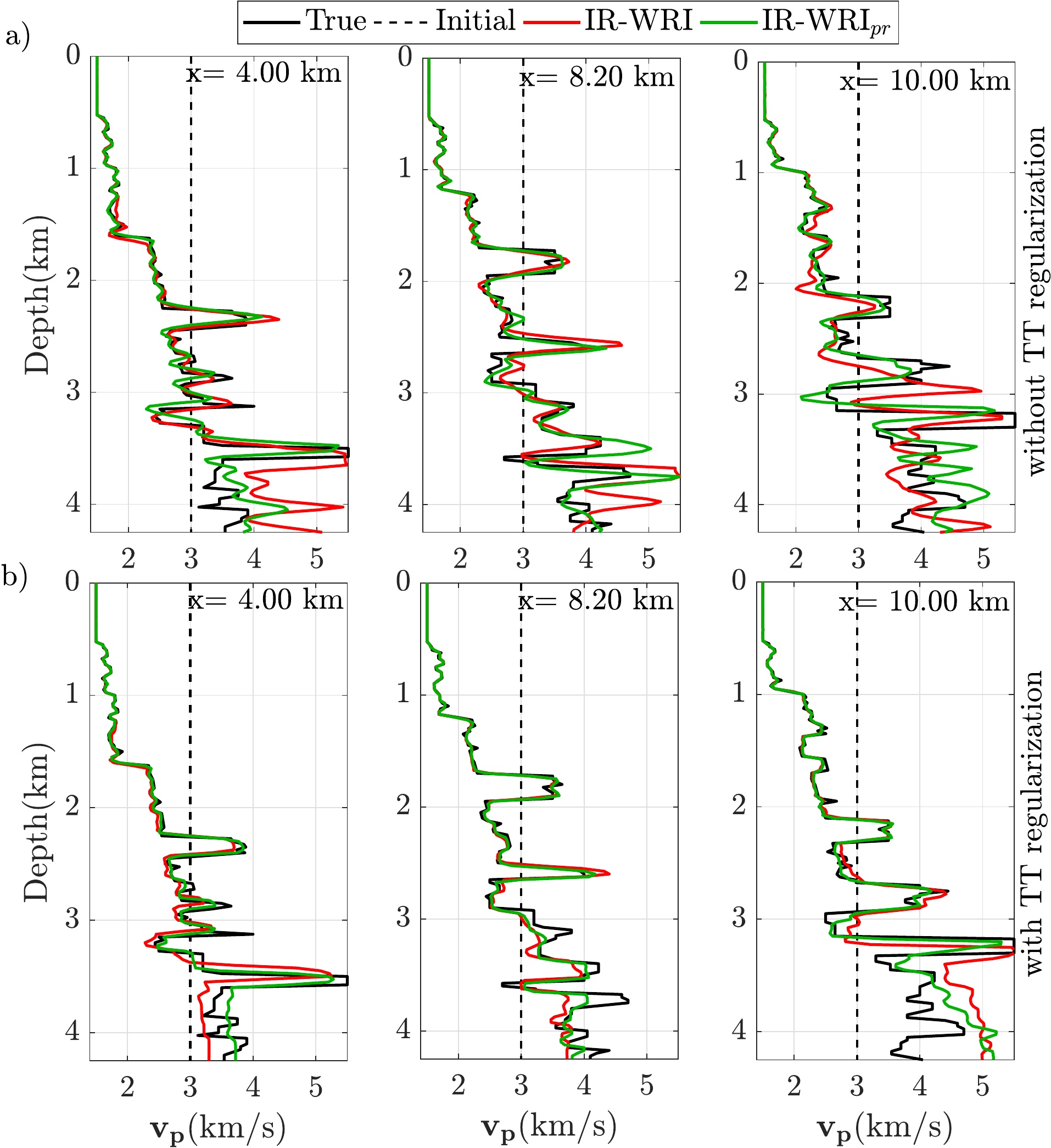}
\caption{Marmousi II example. Direct comparison at $x=4, 8.2$, and $10~km$  between the true velocity model (black), the initial model (dashed line), and the models estimated by IR-WRI (red) and IR-WRI$_{pr}$ (green) (Fig.~\ref{fig:Mar_final}). (a) Without TT regularization. (b) With TT regularization.}
\label{fig:mar_log_final}
\end{figure}

We continue with IR-WRI at higher frequencies using the final IR-WRI and WIPR models of the [3, 3.5]~Hz inversion (Fig.~\ref{fig:mar_model_first}) as initial models. 
The stopping criterion of iterations is a maximum of 30 iterations per batch or 
\begin{equation} \label{Stop}
(\| \bold{A(m}^k)\bold{u}^k-\bold{b}\|_{F} \leq 10^{-3}  \hspace{0.5 cm} \text{and}  \hspace{0.5 cm} \| \bold{Pu}^k-\bold{d}\|_{F} \leq 10^{-5}), 
\end{equation}
where $F$ refers to the Frobenius norm. 
The final IR-WRI and IR-WRI$_{pr}$ velocity models are shown in Fig.~\ref{fig:Mar_final}, while a direct comparison between them and the true model are plotted in Fig.~\ref{fig:mar_log_final} at distances 4, 8.2 and 10~km. IR-WRI performed a total of 382 and 369 iterations without and with TT regularization, respectively (Fig.~\ref{fig:Mar_final}a,c), while IR-WRI$_{pr}$ performed a total of 201 and 234 iterations (Fig.~\ref{fig:Mar_final}b,d). The smaller number of iteration required to satisfy the stopping criterion of iterations performed by IR-WRI$_{pr}$ compared to IR-WRI highlights the role of the initial model built by WIPR during the [3, 3.5]~Hz inversion (Fig.~\ref{fig:mar_model_first}b,d). Without regularization, IR-WRI$_{pr}$ outperforms IR-WRI in the deep part, namely near the reservoir at around 8~km distance and near the salt layers near the ends of the models. Regularization improves significantly both IR-WRI and IR-WRI$_{pr}$ results in these areas. However, IR-WRI$_{pr}$ still provides superior results in the deep part.
%

\subsection{The large contrast 2004 BP salt model}
\begin{figure}[ht!]
\centering
\includegraphics[width=0.48\textwidth]{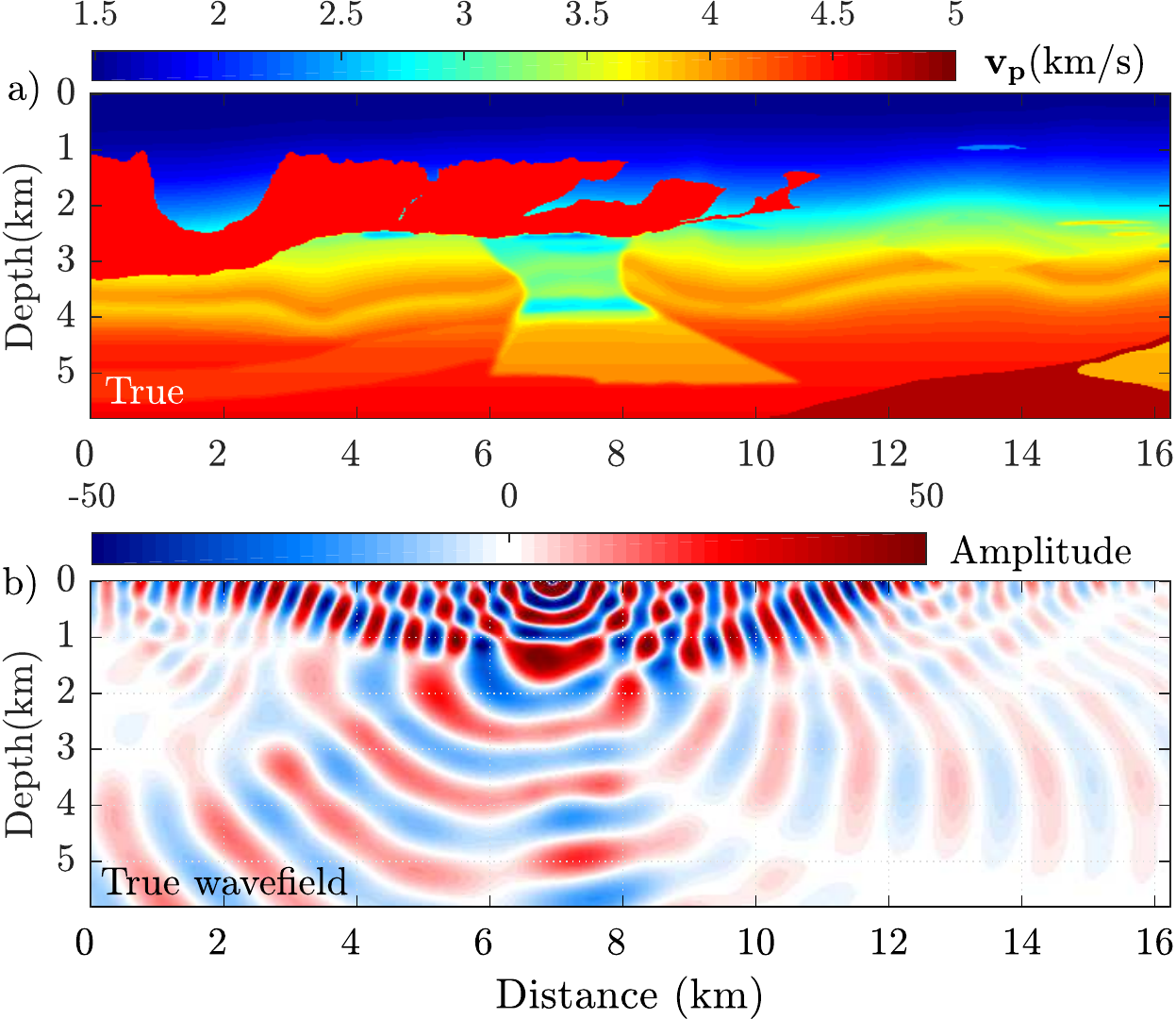}
\caption{(a) True 2004 BP salt velocity model. (b) Real part of the 3-Hz wavefield computed in (a) for the source located at 7 km.}
\label{fig:bp_true}
\end{figure}

We consider now a target of the large contrast 2004 BP salt model located on the left side of the model (Fig. \ref{fig:bp_true}a). The 2004 BP salt model mainly consists of a simple sediment background with a complex rugose multi-valued salt body, sub-salt slow velocity anomalies related to over-pressure zones, and a fast velocity anomaly to the right of the salt body. 
The selected subsurface model is 16250~m wide and 5825~m deep and is discretized with a 25~m grid interval. The surface fixed-spread acquisition consists of 66 sources spaced 250~m apart, and 131 receivers spaced 125~m apart. We perform forward modelling using a 10~Hz Ricker wavelet as a source signature. The starting velocity model for inversion is homogeneous with a velocity of 3 km/s. \\
%
%

\subsubsection{IR-WRI versus WIPR at low frequencies}
\begin{figure}[ht!]
\centering
\includegraphics[width=0.48\textwidth]{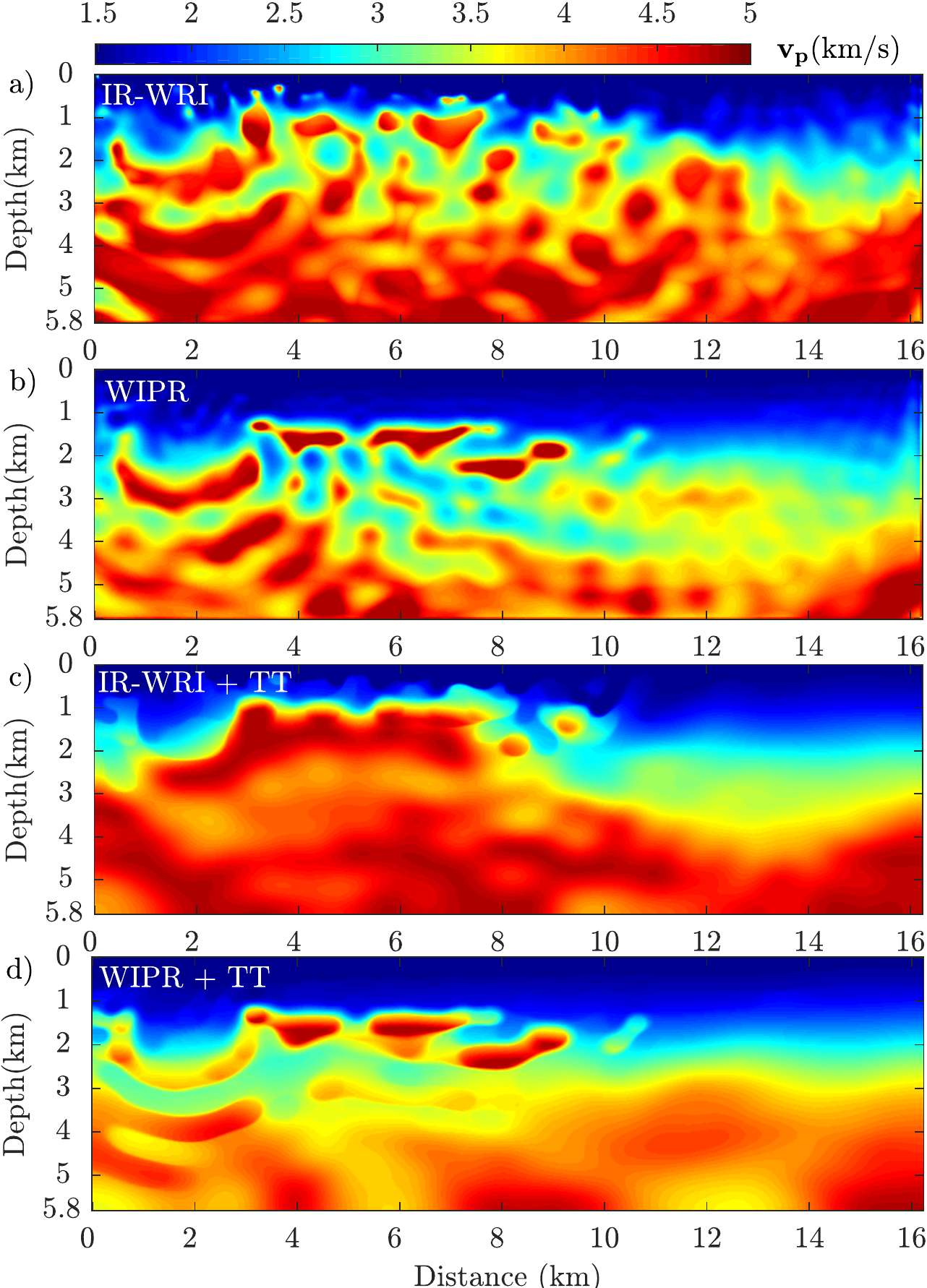}
\caption{BP salt example. Velocity models inferred from the inversion of the \{3,3.5\}~Hz frequencies. (a) IR-WRI, (b) WIPR, (c) IR-WRI with TT regularization, and (d) WIPR with TT regularization. }
\label{fig:bp_model_first}
\end{figure}

\begin{figure}[ht!]
\centering
\includegraphics[width=0.48\textwidth]{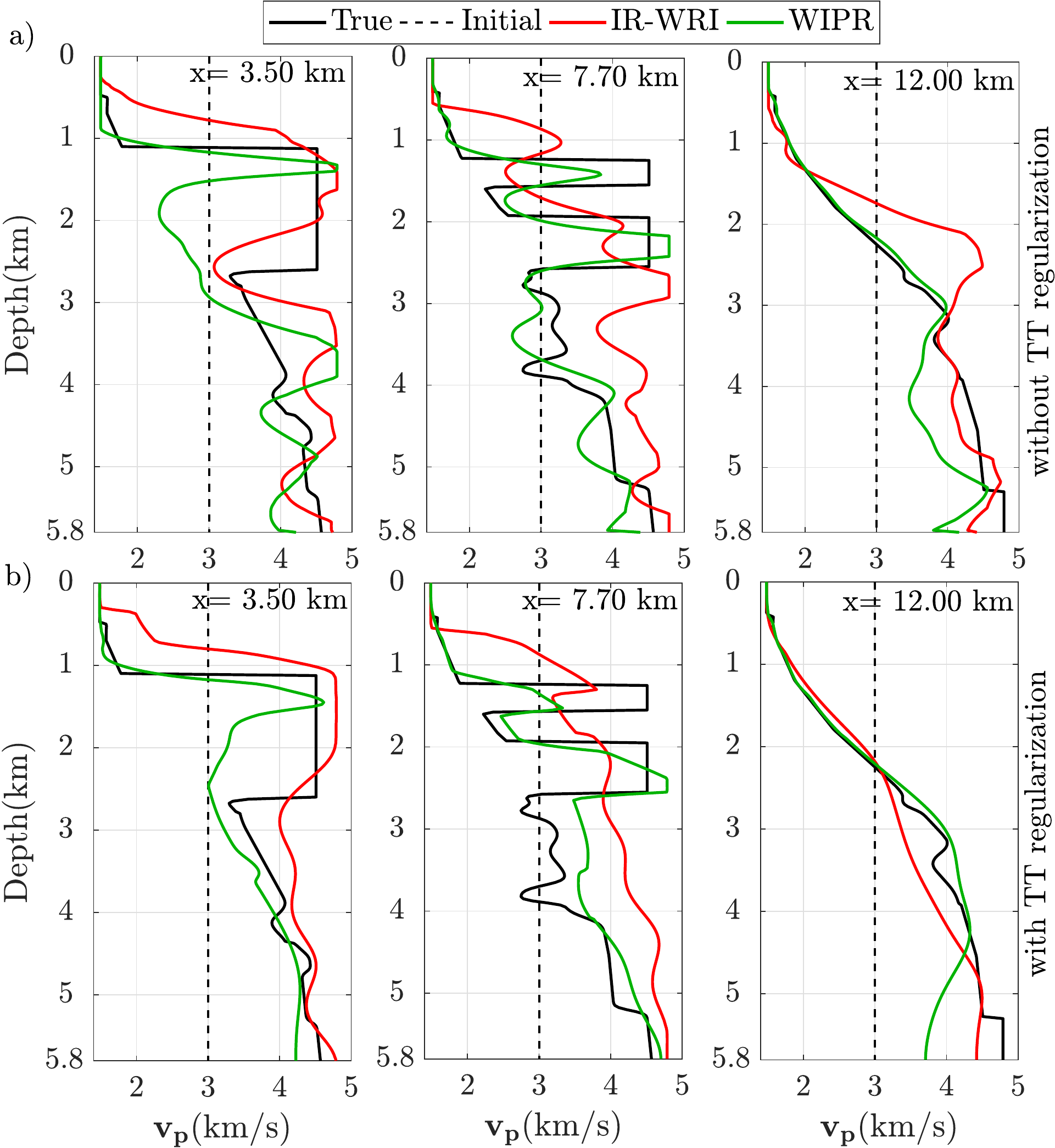}
\caption{BP salt example. Direct comparison along logs located at $x=3.5, 7.7$ and $12~km$ between the true velocity model (black), the initial model (dashed line) and the models estimated by IR-WRI (red) and WIPR (green) shown in Fig.~\ref{fig:bp_model_first}.  (a) Without TT regularization. (b) With TT regularization. }
\label{fig:bp_log_first}
\end{figure}

We first perform 45 iterations of bound-constrained WIPR and IR-WRI for a starting [3,3.5]~Hz frequency batch. Note that we turn on the bound constraints after 20 iterations to discriminate their role from those of the TT regularization and optimization scheme (see later discussion in the text). 
The velocity models reconstructed by bound-constrained IR-WRI and WIPR with and without TT regularization are shown in Fig.~\ref{fig:bp_model_first}, while the direct comparison between the true model, the initial model, and the models reconstructed by IR-WRI and WIPR are shown in Fig.~\ref{fig:bp_log_first} along three vertical logs at horizontal distances of 3.5~km, 7.7~km, and 12~km. 
The most striking feature is that IR-WRI overestimates velocities above the salt and mispositions the top of the salt layers accordingly, while WIPR captures the top salt much more accurately (Fig.~\ref{fig:bp_log_first}a, green versus red curves, x=3.5~km and 7.7~km). This overestimation of the shallow velocities by IR-WRI is also shown in the right part of the model away from the salt (Fig.~\ref{fig:bp_log_first}a, red curve, x=12~km). On the other hand, WIPR fails to reconstruct the bottom part of the thick salt body and overprints unrealistic variations on the subsalt area when TT regularization is not used (Fig.~\ref{fig:bp_model_first}b and Fig.~\ref{fig:bp_log_first}a, green curve, x=3.5~km). These artefacts are likely generated by the increased ill-posedness of WIPR with depth, which itself results from the lack of phase information. This ill-posedness is however, efficiently mitigated by the TT regularization (Fig. \ref{fig:bp_model_first}d). The TV component of the TT regularization helps to sharpen the top and bottom of the salt layers (Fig.~\ref{fig:bp_log_first}b, green curve, x=7.7km), recovers better the geometry of the  salt layer (Fig.~\ref{fig:bp_log_first}b, green curve, x=3.5km), while the Tikhonov component stabilizes the inversion in-depth and drives it toward smooth reconstruction of the poorly-illuminated subsalt structures (Fig.~\ref{fig:bp_log_first}b, green curve, x=12km). During IR-WRI, TT regularization fails to achieve these goals (Fig. \ref{fig:bp_model_first}c) because the mispositioning in depth of the velocity structures generated by the increased nonlinearity of IR-WRI are too heavy to be corrected by the regularization (Fig.~\ref{fig:bp_log_first}b, red curve). 

%
%
\begin{figure}[ht!]
\centering
\includegraphics[width=0.48\textwidth]{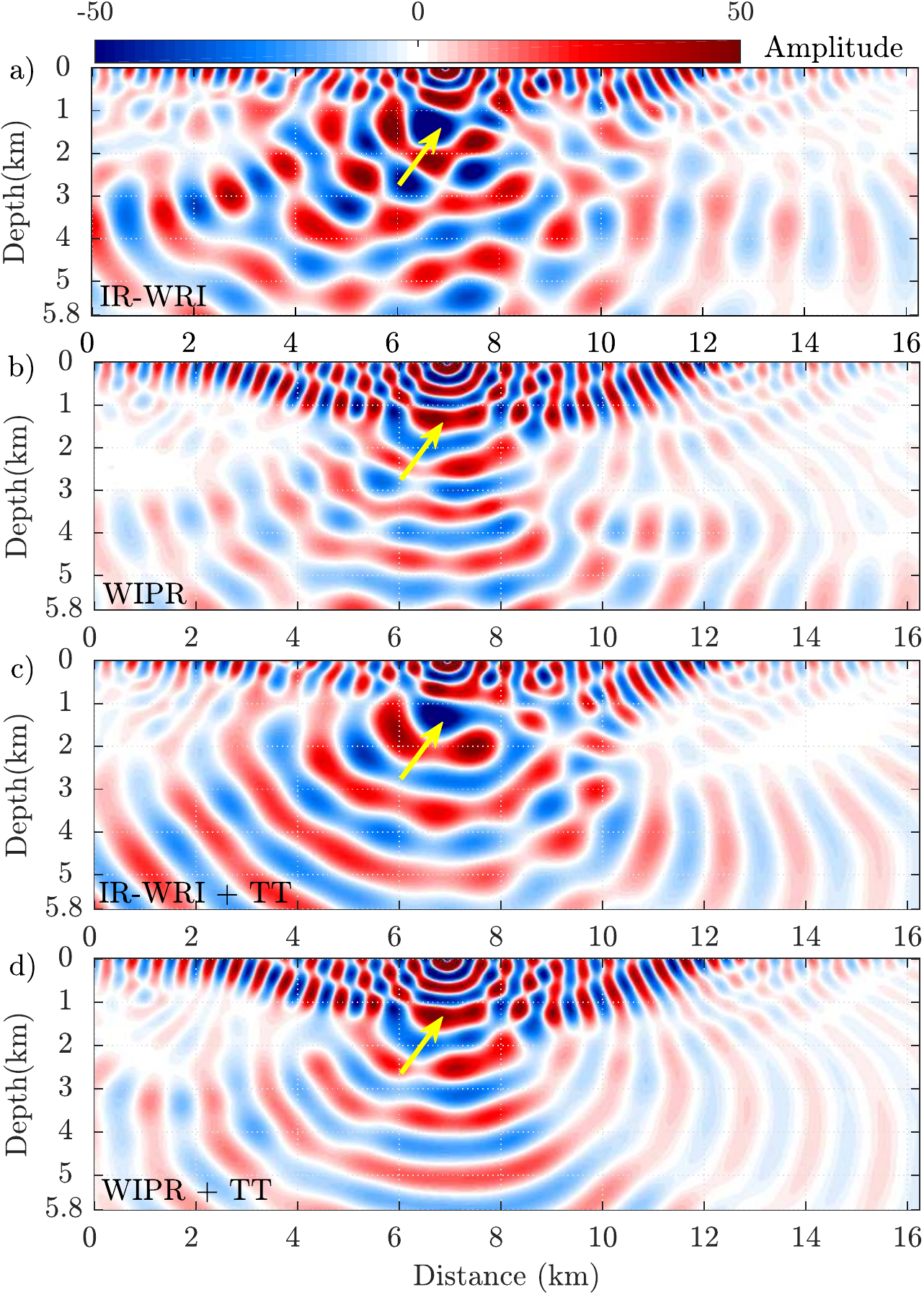}
\caption{BP salt example. Three-Hz wavefields (real part) reconstructed after the first frequency batch inversion.  (a) IR-WRI, (b) WIPR, (c) TT regularized IR-WRI, and (d) TT regularized WIPR. The corresponding velocity models are shown in Fig. \ref{fig:bp_model_first}. The source is located at 7~km distance. These wavefields can be compared with the true wavefield shown in Fig.~\ref{fig:bp_true}b. The arrows point where the incident wavefields interact with the top salt.}
\label{fig:bp_wavefield_first}
\end{figure}

We compare the true 3-Hz wavefield shown in Fig. \ref{fig:bp_true}b with the wavefields reconstructed by IR-WRI and WIPR without and with TT regularization (Fig.~\ref{fig:bp_wavefield_first}). Like the Marmousi test, both IR-WRI and WIPR wavefields are in phase with the true wavefield along the receiver line since we tune the penalty parameter such that the data are matched from the first iteration. However, the WIPR wavefields match much better the true wavefield than the IR-WRI counterpart in the first two kilometers of the model. This is pointed by the arrows in  Fig. \ref{fig:bp_wavefield_first}, which show that the IR-WRI wavefield transmitted into the salt body is out-of-phase with respect to the true wavefield, unlike the WIPR wavefield.\\

%
\begin{figure}[ht!]
\centering
\includegraphics[width=0.48\textwidth]{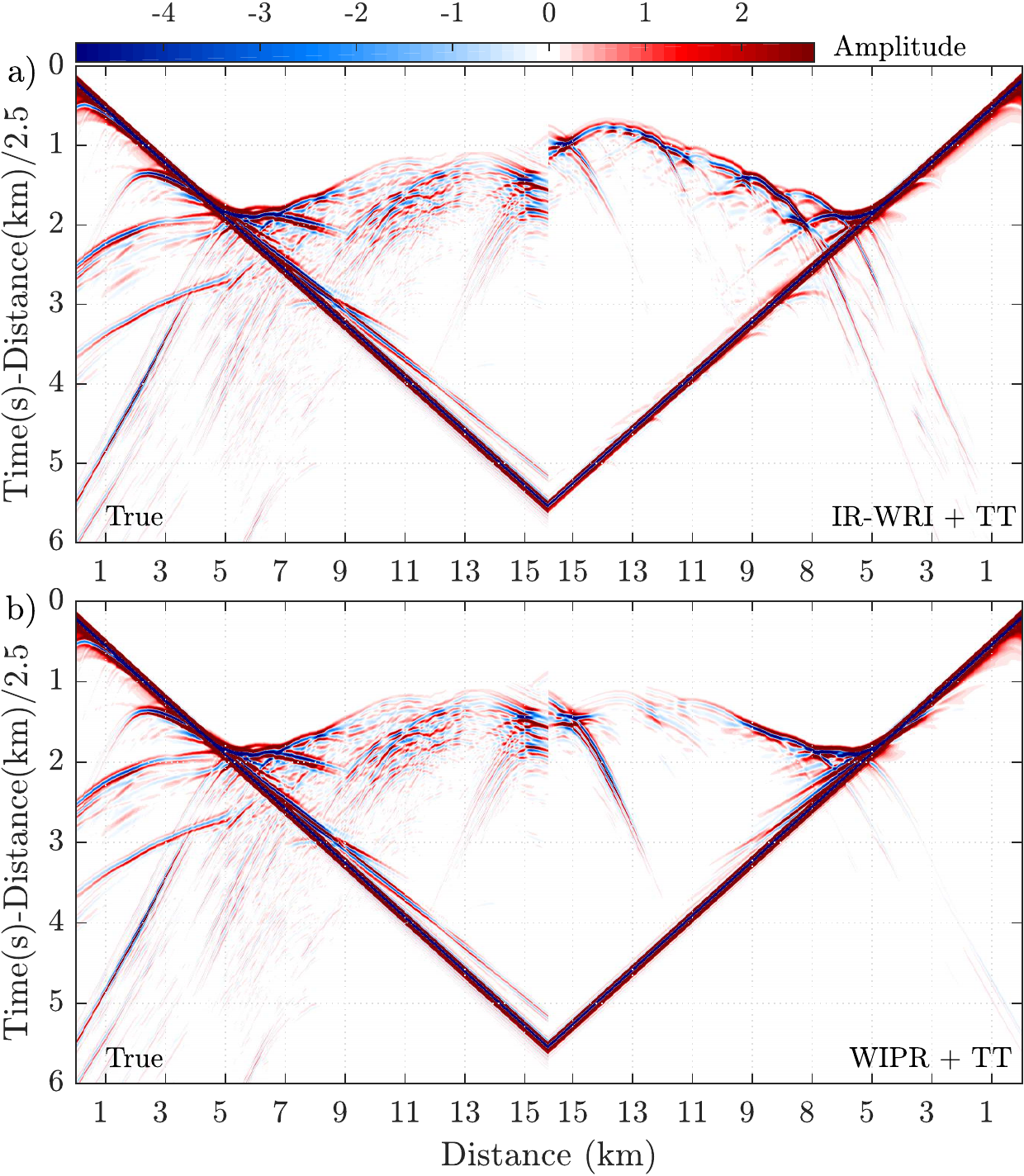}
\caption{BP salt example. Comparison between time-domain seismograms computed in the true velocity model (left panel) and those computed in the velocity models obtained by (a) IR-WRI (Fig.~\ref{fig:bp_model_first}c) and (b) WIPR (Fig.~\ref{fig:bp_model_first}d) (right panel). The seismograms are plotted with a reduction velocity of 2.5 km/s. Note the significant travel-time mismatch in (a) at long offset.}
\label{fig:bp_seismo}
\end{figure}

\vspace{2cm}
We also computed time-domain seismograms in the TT-regularized IR-WRI and WIPR models (Fig.~\ref{fig:bp_model_first}c,d) with a 10-Hz Ricker wavelet for a source located at 15.8~km (Fig. \ref{fig:bp_seismo}). 
We show quite significant travel-time mismatches at long offsets between the true seismograms and those computed in the IR-WRI model (Fig. \ref{fig:bp_seismo}a), which result from the overestimated velocities above the salt shown in Fig. \ref{fig:bp_model_first}c, while the match between the true seismograms and those computed in the WIPR model is far better. This poor data fit achieved by IR-WRI highlights both the inaccuracies of the reconstructed velocity model and the limited accuracy with which the wave equation has been satisfied at the convergence point. \\

%
%
%
%

\begin{figure}[ht!]
\centering
\includegraphics[width=0.48\textwidth]{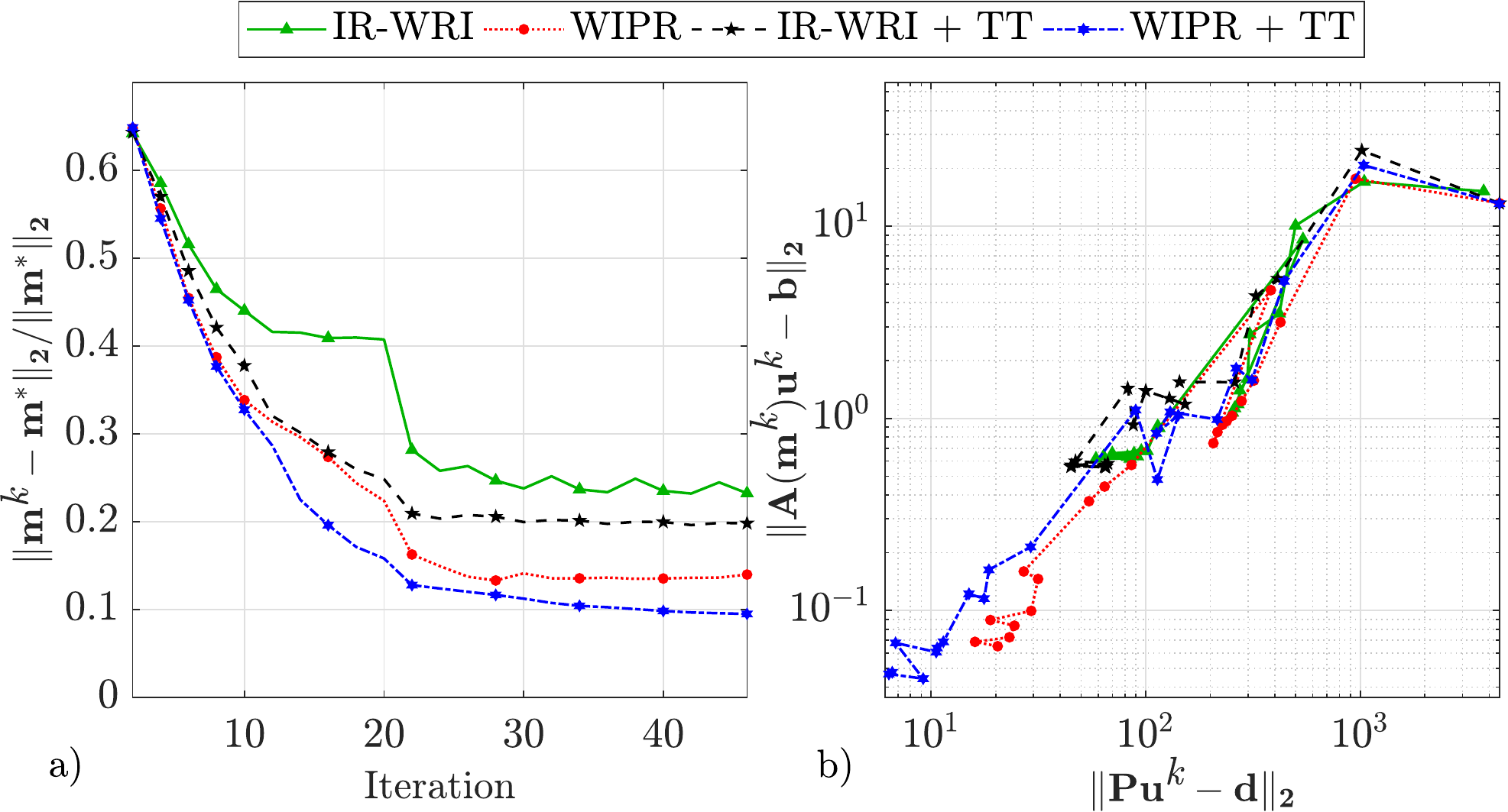}\\
\includegraphics[width=0.48\textwidth]{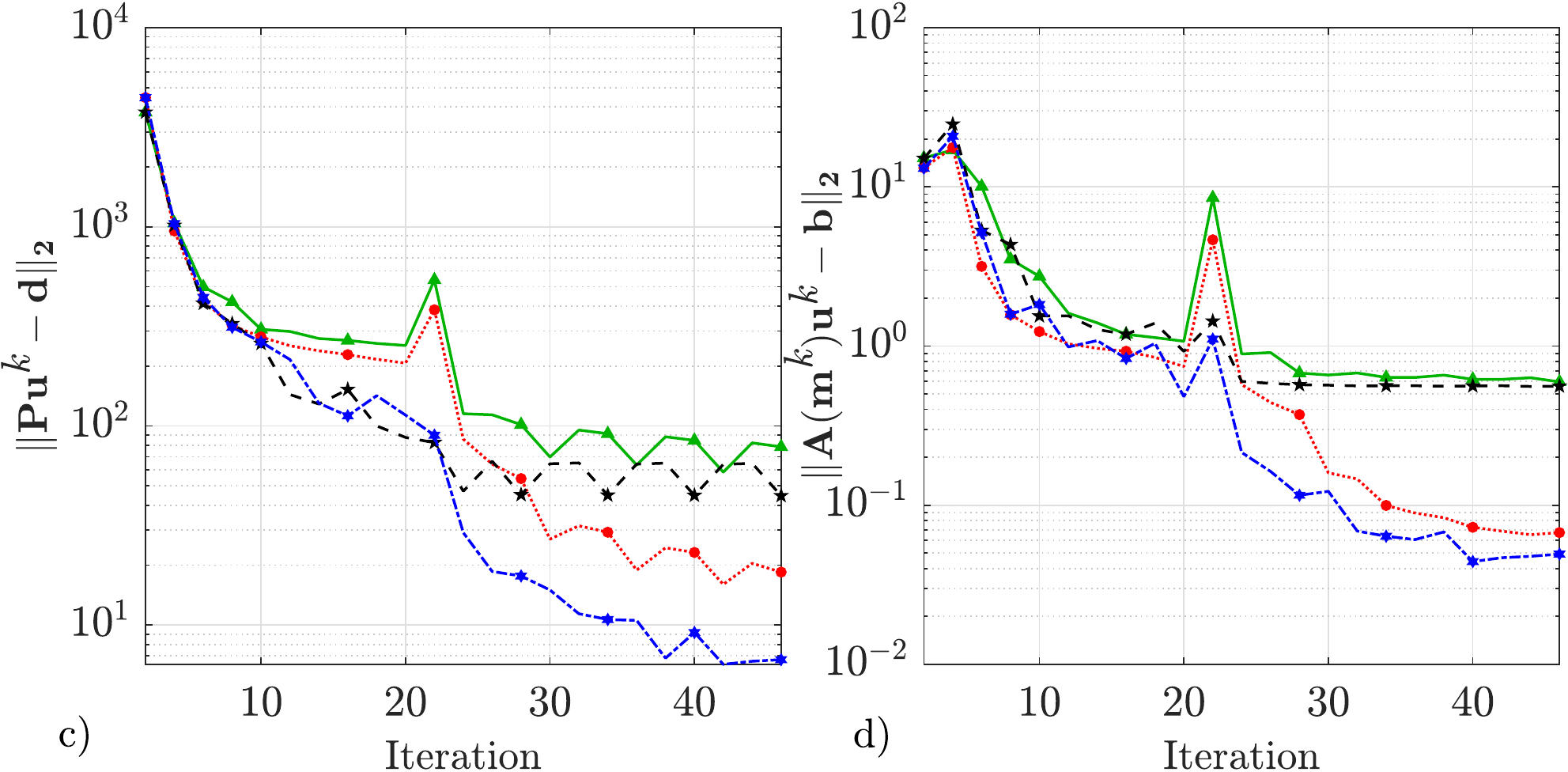}\\
\caption{BP salt example. Convergence path of \{3,3.5\}~Hz inversion. 
(a) Evolution of $\|\bold{m}^k-\bold{m}^*\|_2/\|\bold{m}^*\|_2$ during iterations ($\bold{m}^*$ denotes the true model). (b) Convergence history of the algorithm in the ($\|\bold{Pu}^k-\bold{d}\|_2 - \|\bold{A(m}^k)\bold{u}^k-\bold{b}\|_2$) plane. (c-d) Evolution of (c) $\|\bold{Pu}^k-\bold{d}\|_2$, (d) $\|\bold{A(m}^k)\bold{u}^k-\bold{b}\|_2$ during iterations for IR-WRI (green), WIPR (red), TT regularized IR-WRI (black), and TT regularized WIPR (blue).}
\label{fig:bp_start_res}
\end{figure}

This is further supported by Fig.~\ref{fig:bp_start_res}, which shows the relative model error (Fig.~\ref{fig:bp_start_res}a) and the joint evolution of the data misfit ($\|\bold{Pu}^k-\bold{d}\|_2$) and wave-equation error ($\|\bold{A(m}^k)\bold{u}^k-\bold{b}\|_2$) in iterations (Fig.~\ref{fig:bp_start_res}b). Also, the evolution of the data misfit and wave-equation error in iterations are shown separately in Figs.~\ref{fig:bp_start_res}c and \ref{fig:bp_start_res}d, respectively. First, these error curves show that WIPR converges toward a more accurate velocity model than IR-WRI (Fig.~\ref{fig:bp_start_res}a). Accordingly, it shows that WIPR satisfies more accurately the wave-equation constraint than IR-WRI (Fig.~\ref{fig:bp_start_res}d), although the parameter estimation has been recast as a phase retrieval problem rather than as a wave-equation error minimization. Fig.~\ref{fig:bp_start_res}a also shows that bound constraints have the least impact on the TT regularized WIPR (blue curve) since the descent direction is rather continuous around iteration 20 (when bound constraints are activated) . This highlights the more stable descent direction followed by TT regularized WIPR relative to IR-WRI. \\

%
\subsubsection{IR-WRI versus $\text{IR-WRI}_\textit{pr}$} 
\begin{figure}[ht!]
\centering
\includegraphics[width=0.48\textwidth]{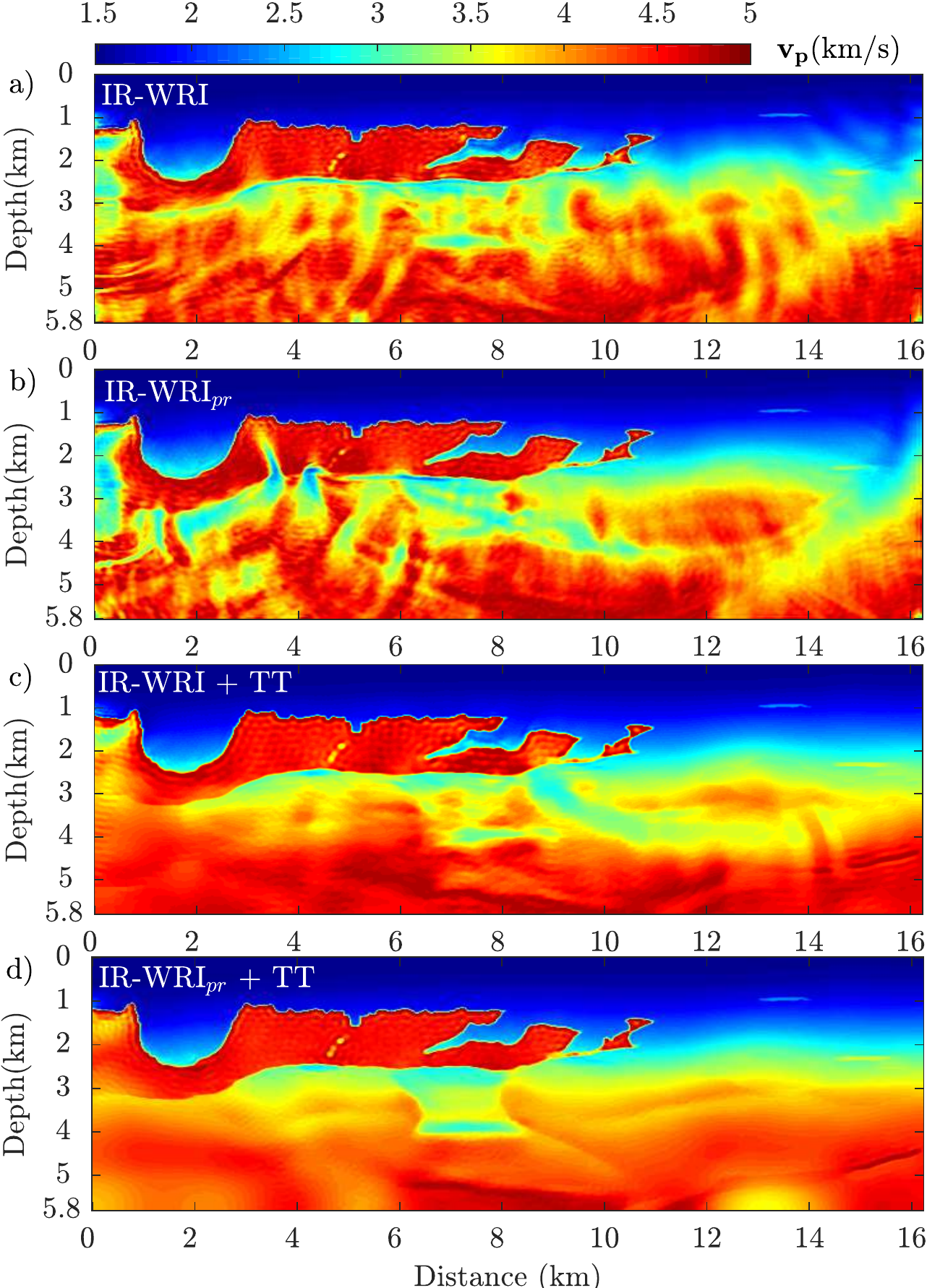}
\caption{BP salt example. Final velocity models estimated by IR-WRI (a,c) and IR-WRI$_{pr}$ (b,d). (a-b) Without TT regularization. (c-d) With TT regularization.}
\label{fig:Bp_final}
\end{figure}

\begin{figure}[ht!]
\centering
\includegraphics[width=0.48\textwidth]{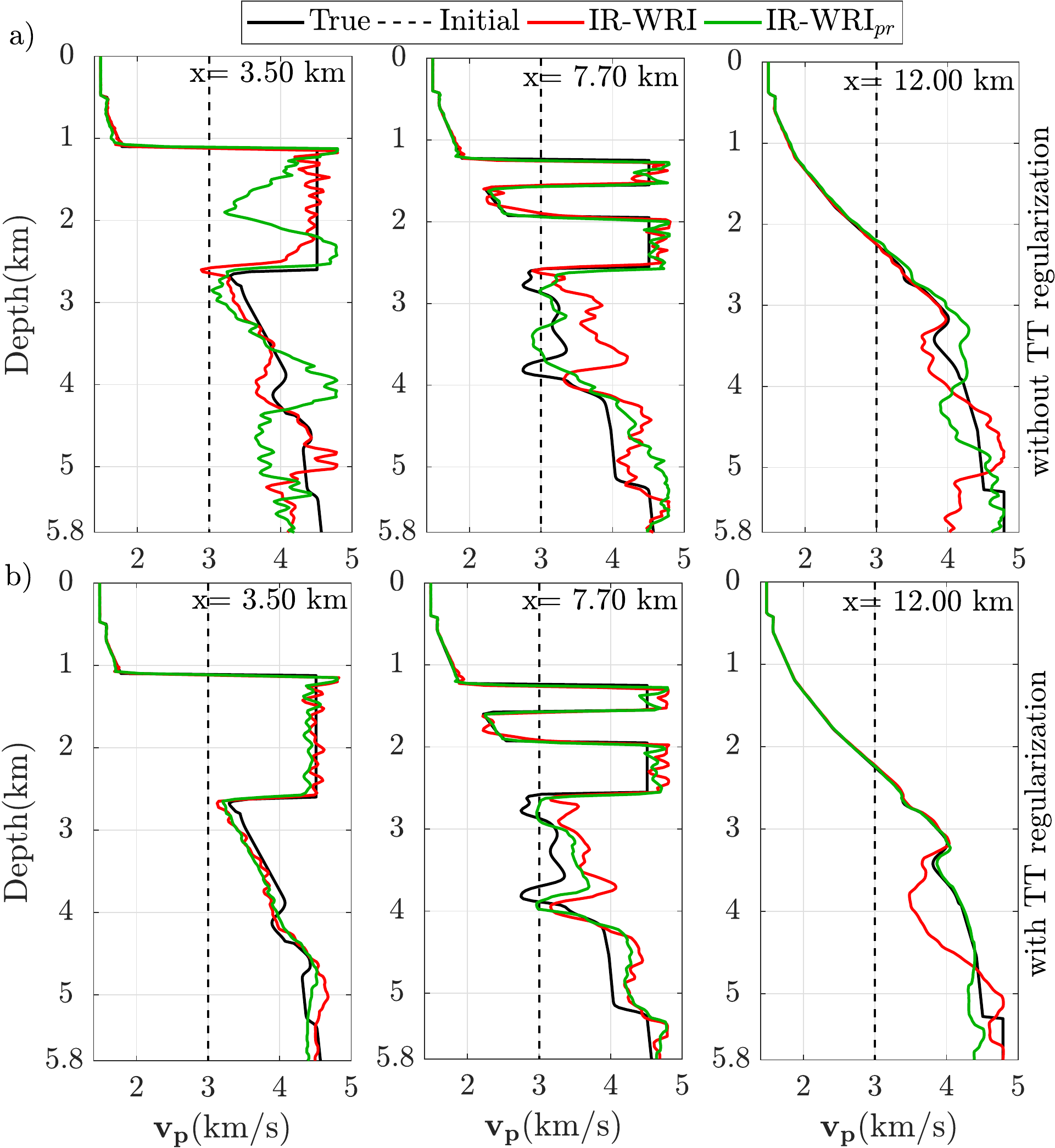}
\caption{BP salt example. Direct comparison along logs located at $x=3.5, 7.7$ and $12~km$ between the true velocity model (black), the initial model (dashed line) and the models estimated by IR-WRI (red) and IR-WRI$_{pr}$ (green) shown in Fig.~\ref{fig:Bp_final}. (a) Without TT. (b) With TT regularization.}
\label{fig:bp_log_final}
\end{figure}

\begin{figure}[ht!]
\centering
\includegraphics[width=0.48\textwidth]{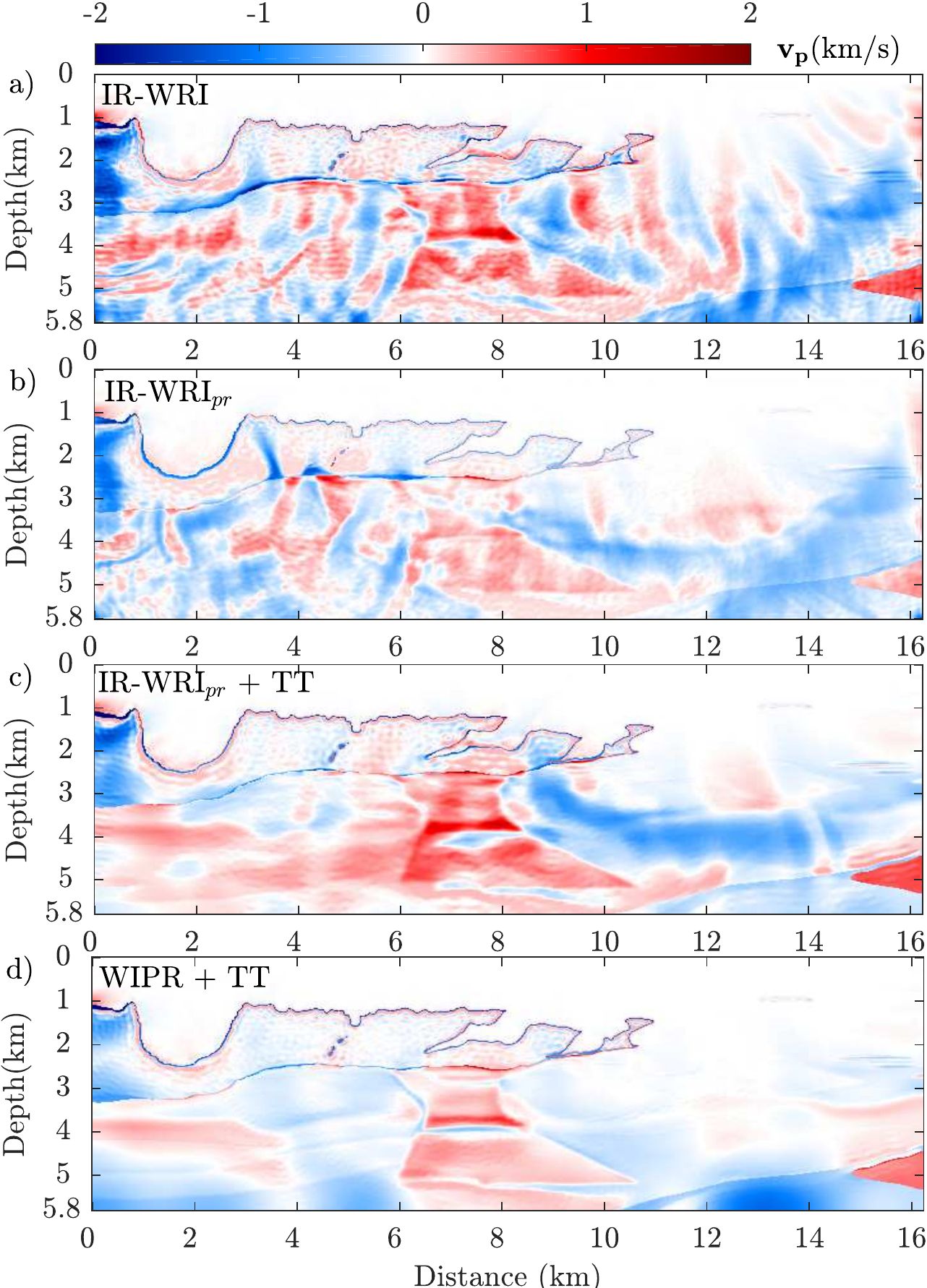}
\caption{BP salt example. Difference between the true velocity model and the velocity models shown in Fig. \ref{fig:Bp_final}. (a,c) IR-WRI. (b,d) IR-WRI$_{pr}$. (a-b) Without TT regularization. (c-d) With TT regularization.}
\label{fig:bp_diff_final}
\end{figure}

We continue with IR-WRI at higher frequencies using the final IR-WRI and WIPR models of the [3, 3.5]~Hz inversion (Fig.~\ref{fig:bp_model_first}) as initial models.
The stopping criterion of iterations is 10 iterations per batch. The final velocity models built by IR-WRI and IR-WRI$_{pr}$ are shown in Fig.~\ref{fig:Bp_final}, while a direct comparison between them and the true model are plotted in Fig. \ref{fig:bp_log_final} at distances 3.5~km, 7.7~km and 12~km. \\
Let's first focus on the results obtained without TT regularization (Fig. \ref{fig:Bp_final}a,b). IR-WRI apparently reconstructs the salt body with good accuracy (Fig. \ref{fig:Bp_final}a), although the shallow kinematic inaccuracies of the starting IR-WRI model shown in Figs. \ref{fig:bp_model_first}a and \ref{fig:bp_log_first}, red curves. This highlights the potential of the search-space expansion implemented in IR-WRI to manage cycle skipping.
The final IR-WRI$_{pr}$ velocity model shows low-velocity artefacts in the salt body at around 3~km distance (Fig. \ref{fig:Bp_final}b), which are probably inherited from the underestimated velocities below the top of the salt shown in Figs.~\ref{fig:bp_model_first}b and \ref{fig:bp_log_first}a, green curve, x=3.5~km. For both inversions, the reconstruction of the subsalt structures are quite noisy. To gain more quantitative insights on the relative accuracy of the IR-WRI and IR-WRI$_{pr}$ models, we show the difference between these velocity models and the true model in Fig.~\ref{fig:bp_diff_final}(a-b) and outline the model error (ME) for the velocity models as well as the least-squares data and source misfit at the convergence point in Table~\ref{taberror}. 
Figure \ref{fig:bp_diff_final}(a-b) clearly shows that the IR-WRI$_{pr}$ model is as a whole more accurate than the IR-WRI counterpart, although the above-mentioned local artefacts.

\begin{table}[ht!]
\begin{tabular}{lcc|c|c|c|}
\cline{4-6}
                                                 &                                                 & \multicolumn{1}{l|}{} & ME  & $\|\bf{A(m)u-b}\|_2$ & $\|\bf{Pu-d}\|_2$ \\ \hline
\multicolumn{1}{|r|}{\multirow{4}{*}{\hspace{-0.2cm}\rotatebox{90}{Iteration}}} \hspace{-0.0cm}& \multicolumn{1}{c|}{\multirow{4}{*}{\hspace{-0.3cm}\rotatebox{90}{No.45}}}     & IR-WRI              & 23.23 & 0.5954               & 7.8550            \\ \cline{3-6} 
\multicolumn{1}{|r|}{}                           & \multicolumn{1}{c|}{}                           & WIPR                      & 13.97 & 0.0672         & 1.8450            \\ \cline{3-6} 
\multicolumn{1}{|r|}{}                           & \multicolumn{1}{c|}{}                           & IR-WRI+TT               & 19.81 & 0.5640          & 4.4630            \\ \cline{3-6} 
\multicolumn{1}{|r|}{}                           & \multicolumn{1}{c|}{}                           & WIPR+TT                 & 9.49  & 0.04926        & 0.6706            \\ \hline
\multicolumn{1}{|l|}{\multirow{4}{*}{\hspace{-0.2cm}\rotatebox{90}{Final}}}     & \multicolumn{1}{c|}{\multirow{4}{*}{\hspace{-0.3cm}\rotatebox{90}{iteration}}} & IR-WRI             & 7.57  & 0.2046               & 2.8911            \\ \cline{3-6} 
\multicolumn{1}{|l|}{}                           & \multicolumn{1}{c|}{}                           & $\text{IR-WRI}_{pr}$      & 6.50  & 0.0647        & 1.6009            \\ \cline{3-6} 
\multicolumn{1}{|l|}{}                           & \multicolumn{1}{c|}{}                           & IR-WRI+TT               & 4.89  & 0.1272          & 1.6470            \\ \cline{3-6} 
\multicolumn{1}{|l|}{}                           & \multicolumn{1}{c|}{}                           & $\text{IR-WRI}_{pr}$+TT & 3.16  & 0.0359          & 0.6350            \\ \hline
\end{tabular}
\caption{Model error (ME) for the velocity models and least-squares data and source misfit at the end of the first frequency batch inversion and at the convergence point of the four inversion tests (with or without phase retrieval, with or without TT regularization). Here, ME is defined as 
$100  \frac{\|\bold{m}_k-\bold{m}_*\|_1}{\|\bold{m}_*\|_1} $, where $\bold{m}_*$ denotes the true model.}
\label{taberror}
\end{table}

The TT regularization removes the local artefacts shown in the salt body of the IR-WRI$_{pr}$ model and improves the subsalt imaging (Fig.~\ref{fig:Bp_final}d) more significantly than for the IR-WRI model (Fig.~\ref{fig:Bp_final}c). The difference between the true model and the IR-WRI/IR-WRI$_{pr}$ models (Fig.~\ref{fig:bp_diff_final}(c-d)) and the model error (ME) provided in Table~\ref{taberror} outline the superior accuracy of the IR-WRI$_{pr}$ model relative to the IR-WRI counterpart. These relative model accuracies are consistent with the data and source misfit achieved by IR-WRI$_{pr}$ and IR-WRI at the convergence point (Table~\ref{taberror}).

%

\section{DISCUSSION} 
With the large contrast BP salt model, we have shown that the phase inaccuracies of the reconstructed wavefields prevent accurate estimation of the shallow sediments and accurate positioning in depth of the top of the salt after the first frequency batch inversion. In contrast, WIPR reconstructs more accurately the velocity gradient in the sedimentary cover as well as the sharp velocity contrast on top of the salt. 
This accurate positioning of the top salt is indeed important to speed up the convergence and improve the solution at higher frequencies. 

\begin{figure}
\centering
\includegraphics[width=0.48\textwidth]{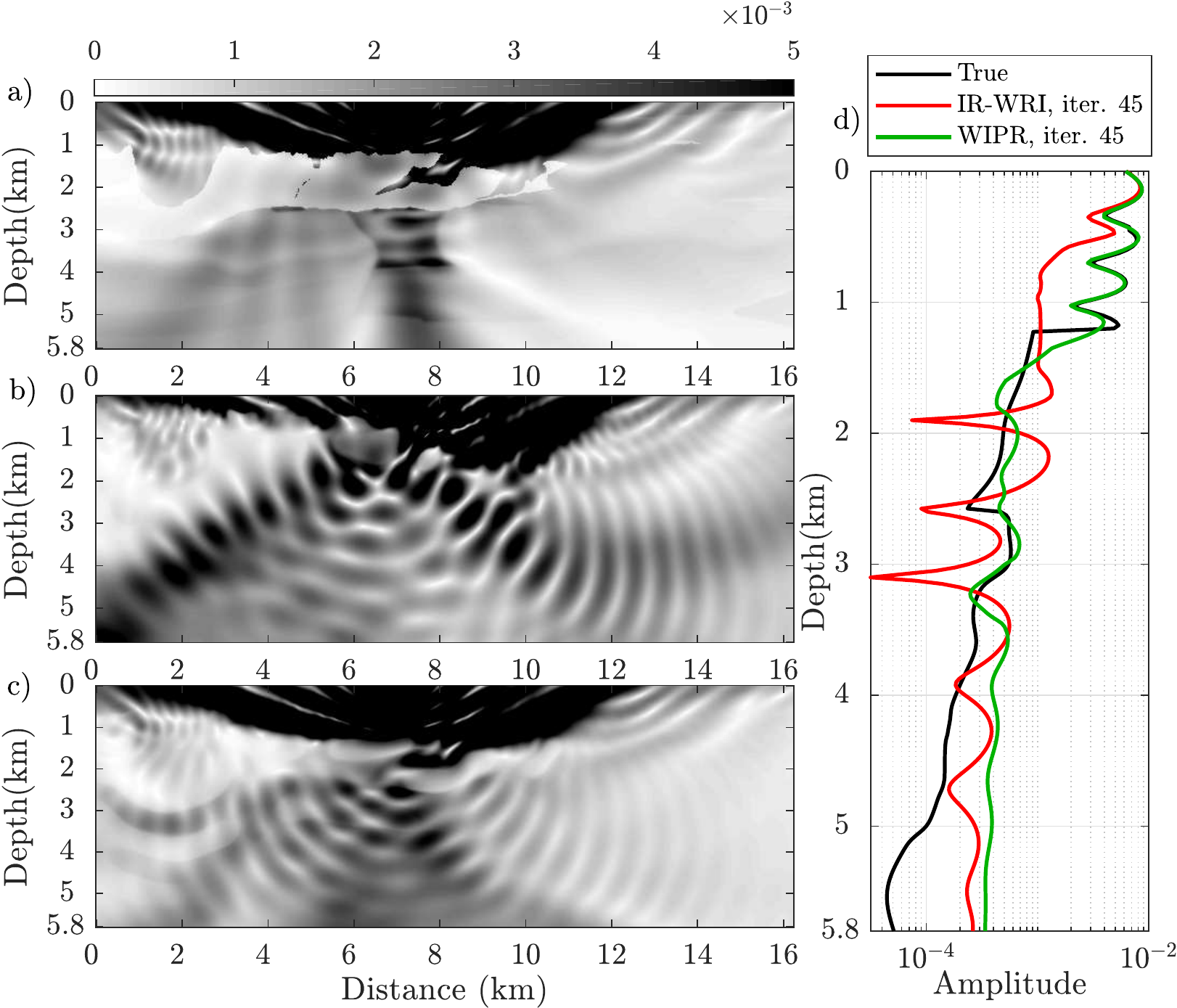}
\caption{Amplitude of $\Delta \bold{u}$ for 2004 BP model. (a) The wavefield $\bold{u}$ is computed in the true velocity model for the 3~Hz frequency. (b-c) The wavefields are computed in the velocity models inferred from TT-regularized (b) IR-WRI (Fig.~\ref{fig:bp_model_first}c) and (c) WIPR (Fig.~\ref{fig:bp_model_first}d)  for the first frequency batch. (d) Comparison between vertical profiles extracted from (a) (black line), (b) (red line), (c) (green line) at a distance of 6.5~km.}
\label{fig:bp_wavefield_amplitude}
\end{figure}

\begin{figure}
\centering
\includegraphics[width=0.48\textwidth]{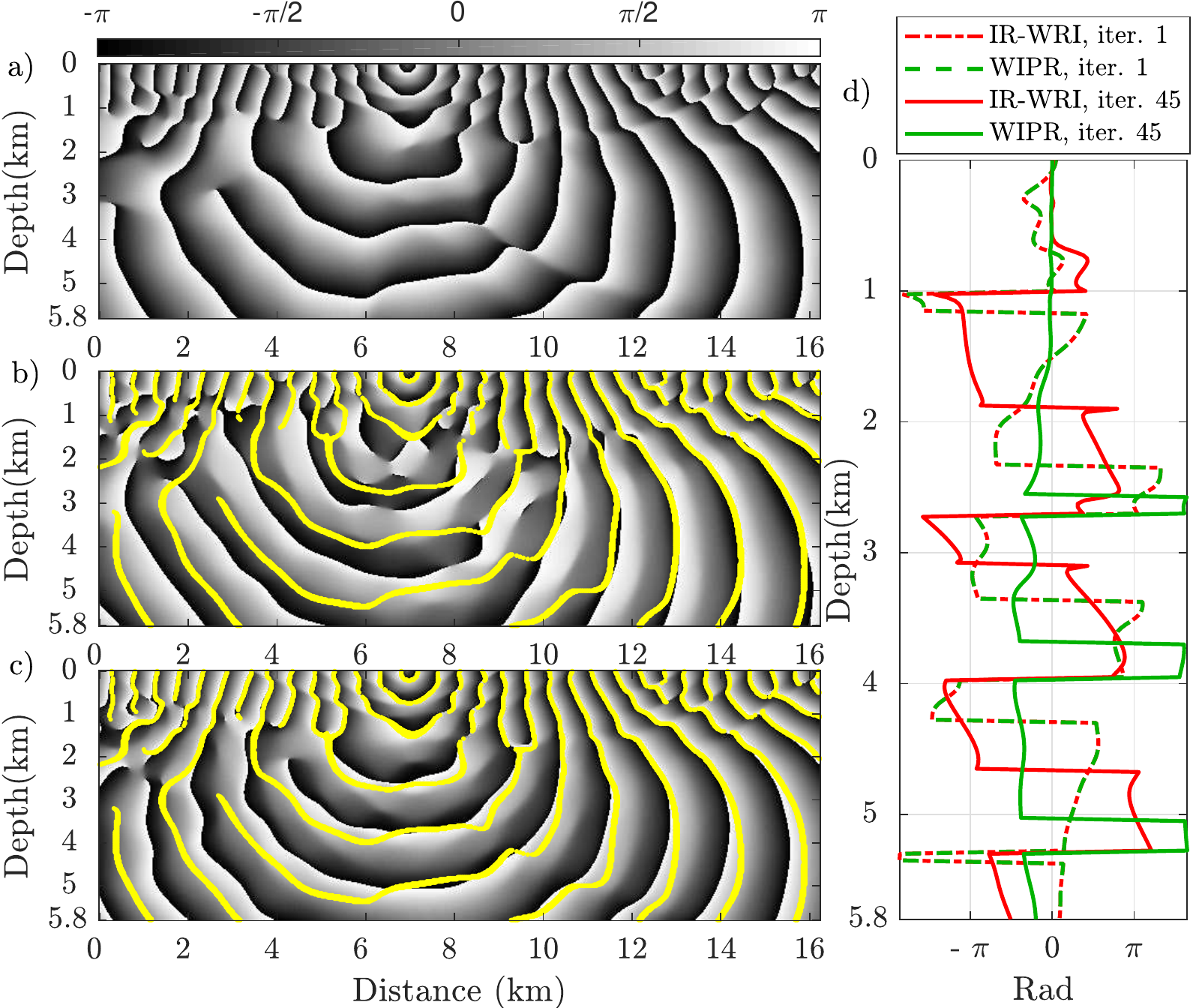}
\caption{Phase of $\Delta \bold{u}$ for 2004 BP model. Same as Fig. \ref{fig:bp_wavefield_amplitude} for the phase of $\Delta \bold{u}$. (d) Difference between the phases computed in the true model and the reconstructed models at the first iteration (dashed red and green lines),  after the first frequency batch inversion for IR-WRI (red line) and WIPR (green line) at a distance of 6.5~km. Note the phase mismatch between the phase computed in the true model and in the IR-WRI model below 0.5~km depth.}
\label{fig:bp_wavefield_phase}
\end{figure}

To further illustrate the sensitivity of IR-WRI to wavefield amplitudes and phases in large contrast media, we show in Figs. \ref{fig:bp_wavefield_amplitude} and \ref{fig:bp_wavefield_phase} the amplitude and phase of $\Delta \bold{u}$ for the 3~Hz frequency and a source located at 7~km of distance. In these figures, $\Delta \bold{u}$ have been computed in the true velocity model, and in the TT regularized IR-WRI and WIPR velocity models obtained after the first frequency batch inversion (Figs. \ref{fig:bp_model_first}c,d). We remind that $\Delta \bold{u}$ is one of the terms forming the right-hand side $\tilde{\bold{y}}$ of the quadratic parameter-estimation subproblem together with the updated source $\bold{b}+{\bold{b}^k}$, Algorithm~\ref{Alg_IR_WIPR}, Line 5. As such, $\Delta \bold{u}$ can be viewed as the observable of the parameter-estimation subproblem. 

The amplitudes of $\Delta \bold{u}$ computed in the true velocity model show the dominant imprint of the reflection and the refraction from the top of the salt (Fig. \ref{fig:bp_wavefield_amplitude}a, shallow black zone). Also, the geometry of the seismically-transparent salt body is almost perfectly delineated. This shows that the amplitudes of $\Delta \bold{u}$ intrinsically embeds high-resolution information on the subsurface. 

This high-resolution potential can be simply illustrated with the wave-equation bilinearity as discussed below. Let us remind the wave equation for a single source and a single frequency
\begin{equation}   \label{AA}
\bold{\Delta}\bold{u} + \omega^2  \text{diag}(\bold{m}) \bold{u}= \bold{b}
\end{equation}
which leads to the following linear equation of the parameter-estimation subproblem: 
\begin{equation}
\underbrace{\omega^2  \text{diag}(\bold{m})\bold{u}}_{\bold{L(u})\bold{m}}=
\underbrace{\bold{b}-\bold{\Delta}\bold{u}}_{\bold{y}(\bold{u})}.
\label{eqwes}
\end{equation}
Compared to equation \ref{A}, we drop the matrix $\bold{B}$ and $\bold{C}$ which are not important for the discussion. 
The equality implies that $|\bold{L(u})\bold{m}|=|\bold{y}(\bold{u})|$ and 
$\angle\bold{L(u})\bold{m}=\angle\bold{y}(\bold{u})$. If $\bold{u}$ is the true wavefield (if it would be recorded everywhere) then we would be able to recover the true model exactly by minimizing
\begin{align} \label{WPR}
\|\bold{L(u})\bold{m} - \bold{y}(\bold{u})\|_2^2
= \|(|\bold{L(u})\bold{m}| - |\bold{y}(\bold{u})|)e^{j\angle \bold{L(u})\bold{m}}\|_2^2\nonumber \\
= \|\bold{L(u})\bold{m} - |\bold{y}(\bold{u})|e^{j\angle \bold{L(u})\bold{m}}\|_2^2.
\end{align}
It is seen that, basically, the true model is the solution of a weighted phase retrieval problem (as indicated by the middle term) where the complex exponential of the phase serves as the weight.  

Indeed, for an approximate wavefield which is obtained from eq. \eqref{u} with a rough velocity model or a large penalty parameter the equality in eq. \eqref{WPR} is not necessarily satisfied because, in this case, $\angle\bold{L(u})\bold{m}\neq \angle\bold{y}(\bold{u})$. 
The proposed WIPR algorithm simply forces this equality to be satisfied at each iteration by aligning the phase of $\bold{y}$ with that of $\bold{L(u})\bold{m}$ (right term in equation~\ref{WPR}). 

The ability to reconstruct the true model from amplitude information only can be illustrated more explicitly if we assume that $\bold{m}$ is real in equation \ref{AA} (non attenuating medium).
The matrix $\bold{L(u})$ is diagonal and contains the so-called virtual source (equation \ref{eqwes}). 
\begin{equation}
\bold{m} = \frac{1}{\omega^2} \left(\frac{1}{\bold{u}}\right) \circ \bold{y}(\bold{u}),
\label{eqbilinpa}
\end{equation}
where $\circ $ denotes element wise multiplication.
Indeed, we can equally reconstruct the true $\bold{m}$ from the magnitude of $\bold{y}$ as
\begin{equation}
\bold{m} = \frac{1}{\omega^2}  \left(\frac{1}{|\bold{u}|}\right) \circ  | \bold{y}(\bold{u}) |,
\label{eqbilina}
\end{equation}
where $1 / \bold{u}$ and $1 / |\bold{u}|$  are element-wise reciprocal of $\bold{u}$ and $|\bold{u}|$, respectively.
The Marmousi II model computed with equations \ref{eqbilinpa} and \ref{eqbilina} by using the true wavefield (Fig.~\ref{fig:mar_true}b) is shown in Fig.~\ref{fig:mar_true}c and matches exactly the true one (Fig.~\ref{fig:mar_true}a).

The Fig.~\ref{fig:bp_wavefield_amplitude}c, as well as the vertical profile extracted at 6.5~km of distance, confirm that WIPR manages to match the amplitudes of the true $\Delta \bold{u}$ down to the top of the salt at 1.25~km depth, while IR-WRI (Fig.~\ref{fig:bp_wavefield_amplitude}b) allows for an amplitude fit down to a maximum depth of 0.5~km only, consistently with the accuracy of the velocity fields reconstructed by IR-WRI and WIPR (Fig. \ref{fig:bp_model_first}c,d). 

It is also instructive to look at the phase of $\Delta \bold{u}$ (Fig.~\ref{fig:bp_wavefield_phase}). To facilitate the comparison of the phases computed in the true model (Fig.~\ref{fig:bp_wavefield_phase}a) and in the reconstructed models (Fig.~\ref{fig:bp_wavefield_phase}b,c), we superimpose in Fig.~\ref{fig:bp_wavefield_phase}b,c the contours of the phase computed in the true model at the spatial positions where the phase wraparounds, namely when it jumps from $\pm \pi$ to  $\mp \pi$  (Fig.~\ref{fig:bp_wavefield_phase}b,c). We clearly show that WIPR reconstructs the phase of $\Delta \bold{u}$ much more accurately than IR-WRI around the salt body. This is further illustrated by the difference between the phases computed in the true model and the reconstructed models along the vertical profile at 6.5~km distance (Fig.~\ref{fig:bp_wavefield_phase}d). The phase difference associated with IR-WRI (solid red curve) reaches a value of $-\pi$ at the top of the salt (1.25~km depth) showing that at these depths the two wavefields are cycle skipped, while the  phase differences associated with WIPR (solid green curve) remains far below this limit down to around 2.7~km depth. 

Conversely, it is worth noting that IR-WRI matches better the phase of the true wavefield than WIPR in the smooth bottom-right part of the model (Compare Figs.~\ref{fig:bp_wavefield_phase}b and \ref{fig:bp_wavefield_phase}c). This highlights that the sensitivity of the inversion to amplitudes decreases rapidly with depth due to lack of illumination, geometrical spreading effects and energy partitioning at interfaces. 
We also want to stress that, although we update the wavespeeds from the amplitudes of $\bold{y}$, we update $\bold{y}$ at each iteration with the phase and amplitude of the wave equation error (the source residuals) in the framework of the ADMM optimization (Algorithm \ref{Alg_IR_WIPR}). This right-hand side updating re-injects at each iteration the phase error as a defect correction term in the optimization. 

\begin{figure}
\centering
\includegraphics[width=0.48\textwidth]{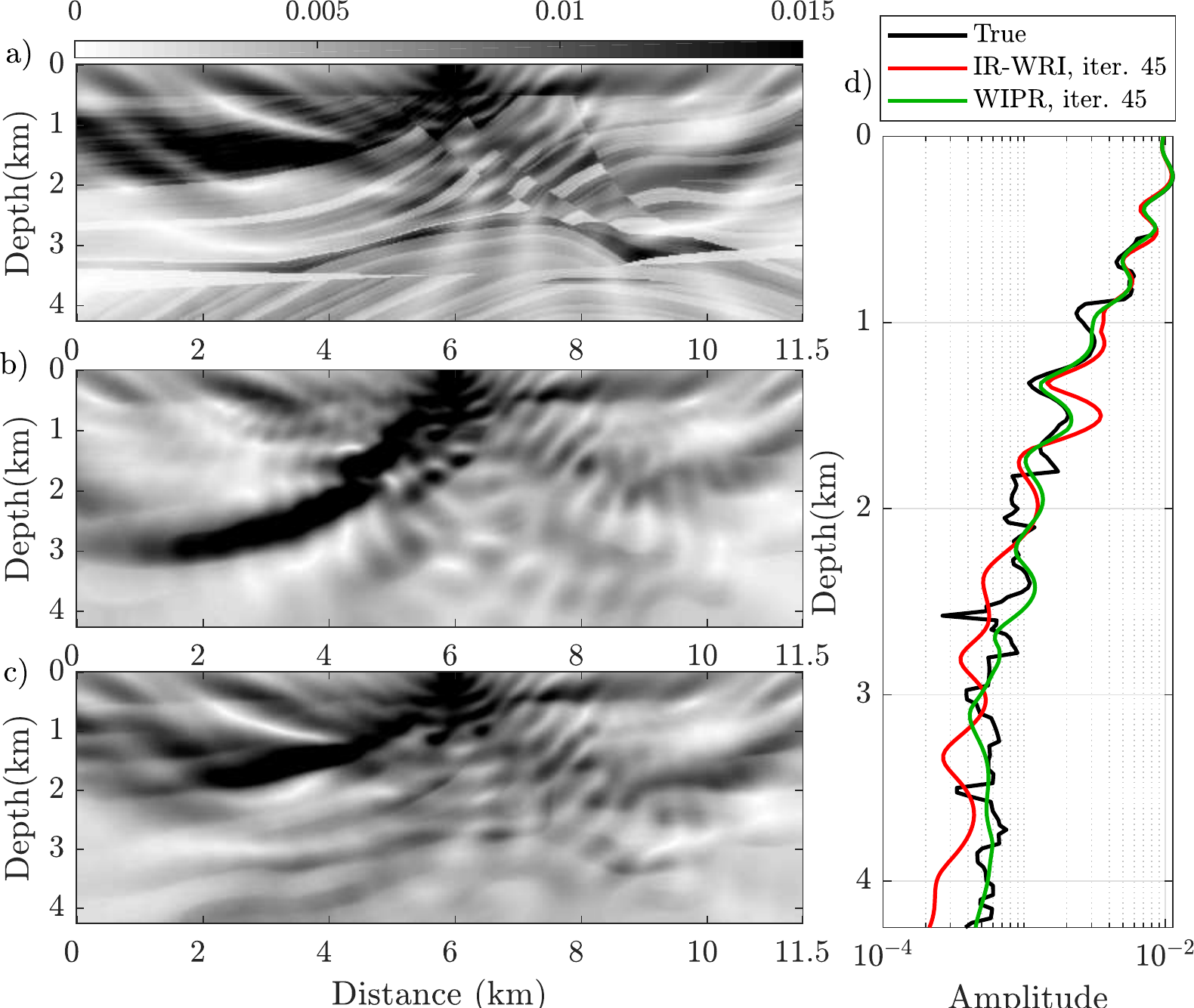}
\caption{Amplitude of $\Delta \bold{u}$ for Marmousi II model. (a) The wavefield $\bold{u}$ is computed in the true velocity model for the 3~Hz frequency. (b-c) The wavefields are computed in the velocity models inferred from TT-regularized (b) IR-WRI (Fig.~\ref{fig:mar_model_first}c) and (c) WIPR (Fig.~\ref{fig:mar_model_first}d)  for the first frequency batch. (d) Comparison between vertical profiles extracted from (a) (black line), (b) (red line), (c) (green line) at a distance of 6.0~km.}
\label{fig:mar_wavefield_amplitude}
\end{figure}

\begin{figure}
\centering
\includegraphics[width=0.48\textwidth]{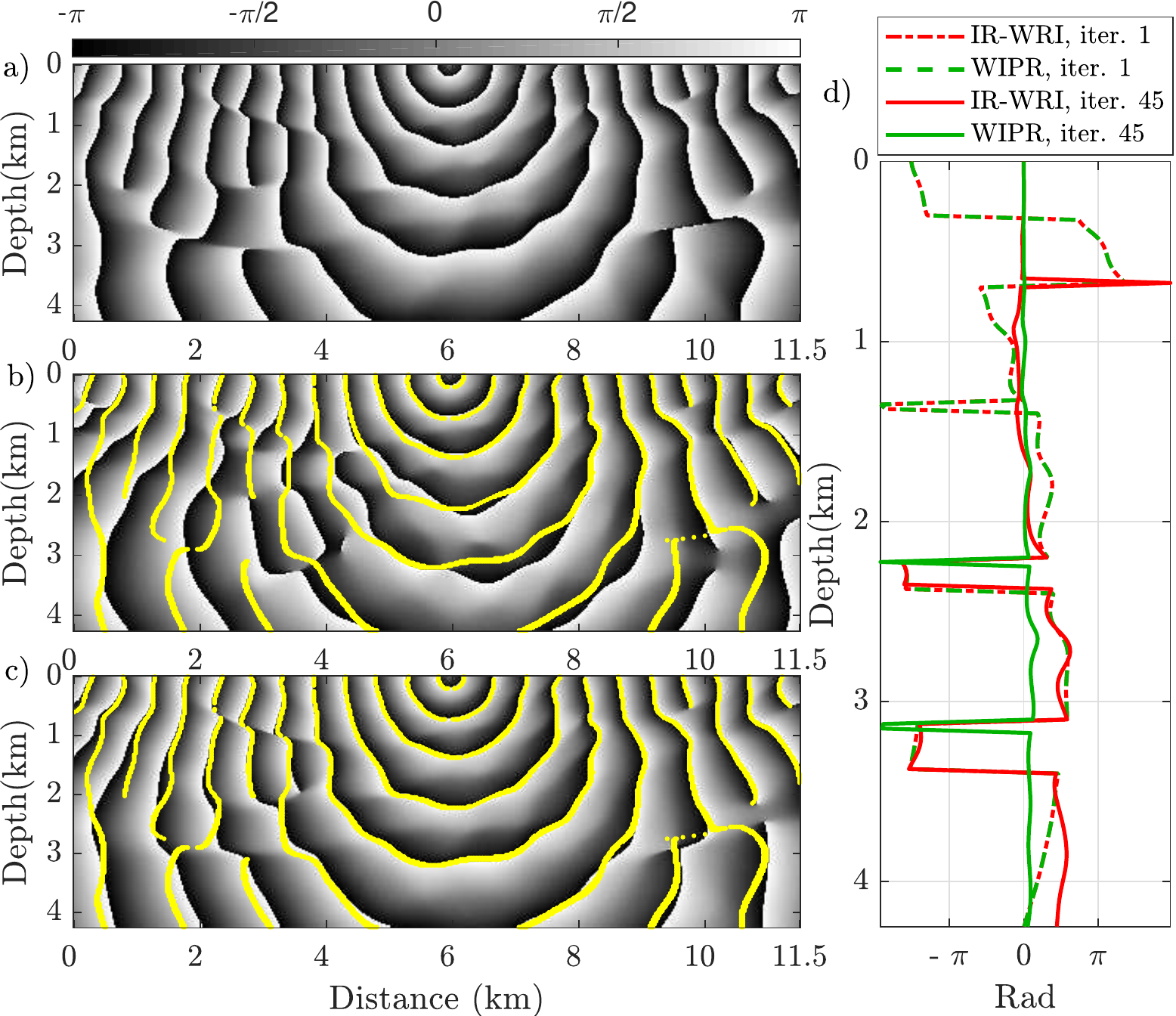}
\caption{Phase of $\Delta \bold{u}$ for Marmousi II model. Same as Fig. \ref{fig:mar_wavefield_amplitude} for the phase of $\Delta \bold{u}$. (d) Difference between the phases computed in the true model and the reconstructed models at the first iteration (dashed red and green lines), after the first frequency batch inversion for IR-WRI (red line) and WIPR (green line) at a distance of 6.0~km.}
\label{fig:mar_wavefield_phase}
\end{figure}

For sake of completeness, we show also the amplitude and phase of $\Delta \bold{u}$ for the Marmousi case study in Figs.~\ref{fig:mar_wavefield_amplitude} and \ref{fig:mar_wavefield_phase}. Fig.~\ref{fig:mar_wavefield_amplitude} shows how WIPR reproduces more accurately the wide-angle reflection from the dipping layer on the left side of the source than IR-WRI. Fig.~\ref{fig:mar_wavefield_phase} shows that WIPR better reconstructs the true phase of $\Delta \bold{u}$ than IR-WRI.

Although the application of the phase retrieval algorithm has been limited to the scalar Helmholtz equation in this study, the method has been formulated such that more complex physics and multi-parameter reconstruction can be viewed (the reader is referred to \citet{Aghamiry_2019_JEO} and \citet{Aghamiry_2019_AMW} for application of IR-WRI in visco-acoustic VTI media). The key feature allowing for these extensions is the linearization of the non-convex phase retrieval problem performed by the MM approach, as reviewed in Appendix \ref{Appa}. 

%
%

%

%

\section{CONCLUSIONS} 
We present a preliminary application of bound-constrained and TT regularized WIPR. Using the large contrast BP salt model, we show that, when the inversion is started from scratch (homogeneous initial velocity model), WIPR improves significantly the reconstruction of the sedimentary background and the top of the salt during the early iterations of the inversion (that is, when the wavefields are not yet accurately reconstructed far away from the receivers). Sparsity-promoting regularization is a necessary ingredient to stabilize the phase retrieval inversion at greater depths where the amplitude information becomes more challenging to extract. Sparsity-promoting regularization combined with phase retrieval maybe also two complementary tools to relax the need for finely-sampled stationary-recording acquisitions. Phase retrieval can also be implemented in classical FWI. Also, WIPR can be theoretically extended to multi-parameter reconstruction in attenuating media. Further theoretical and numerical works are, however, necessary to flesh out the potential and limits of phase retrieval in the FWI technology.  

\section{ACKNOWLEDGEMENTS}  
We would like to thank editors Huajian Yao, Louise Alexander, and reviewers Hejun Zhu and an anonymous for their comments which help improving the manuscript. This study was partially funded by the SEISCOPE consortium (\textit{http://seiscope2.osug.fr}), sponsored by AKERBP, CGG, CHEVRON, EQUINOR, EXXON-MOBIL, JGI, PETROBRAS, SCHLUMBERGER, SHELL, SINOPEC, and TOTAL. This study was granted access to the HPC resources of SIGAMM infrastructure (\textit{https://www.oca.eu/fr/ mesocentre-sigamm}) and CINES/IDRIS/TGCC under the allocation A0050410596 made by GENCI.

\bibliographystyle{seg}

\begin{thebibliography}{48}
\expandafter\ifx\csname natexlab\endcsname\relax\def\natexlab#1{#1}\fi

\bibitem[Aghamiry et~al.(2019{\natexlab{a}})Aghamiry, Gholami, \&
  Operto]{Aghamiry_2019_AMW}
Aghamiry, H., Gholami, A., \& Operto, S., 2019{\natexlab{a}}.
\newblock {ADMM}-based multi-parameter wavefield reconstruction inversion in
  {VTI} acoustic media with {TV} regularization, {\it Geophysical Journal
  International\/}, {\bf 219}(2), 1316--1333.

\bibitem[Aghamiry et~al.(2019{\natexlab{b}})Aghamiry, Gholami, \&
  Operto]{Aghamiry_2019_CRO}
Aghamiry, H., Gholami, A., \& Operto, S., 2019{\natexlab{b}}.
\newblock Compound regularization of {Full-Waveform Inversion} for imaging
  piecewise media, {\it IEEE Transactions on Geoscience and Remote Sensing\/},
  {\bf accepted}.

\bibitem[Aghamiry et~al.(2019{\natexlab{c}})Aghamiry, Gholami, \&
  Operto]{Aghamiry_2019_IBC}
Aghamiry, H., Gholami, A., \& Operto, S., 2019{\natexlab{c}}.
\newblock Implementing bound constraints and total-variation regularization in
  extended full waveform inversion with the alternating direction method of
  multiplier: application to large contrast media, {\it Geophysical Journal
  International\/}, {\bf 218}(2), 855--872.

\bibitem[Aghamiry et~al.(2019{\natexlab{d}})Aghamiry, Gholami, \&
  Operto]{Aghamiry_2019_IWR}
Aghamiry, H., Gholami, A., \& Operto, S., 2019{\natexlab{d}}.
\newblock Improving full-waveform inversion by wavefield reconstruction with
  alternating direction method of multipliers, {\it Geophysics\/}, {\bf 84(1)},
  R139--R162.

\bibitem[Aghamiry et~al.(2019{\natexlab{e}})Aghamiry, Gholami, \&
  Operto]{Aghamiry_2019_JEO}
Aghamiry, H., Gholami, A., \& Operto, S., 2019{\natexlab{e}}.
\newblock Joint estimation of velocity and attenuation by frequency-domain
  {TV}-regularized wavefield reconstruction inversion, in {\em 81$^{th}$ Annual
  EAGE Meeting (London) - WS01: Attenuation: Challenges in Modelling and
  Imaging at the Exploration Scale\/}.

\bibitem[B\'erenger(1994)]{Berenger_1994_PML}
B\'erenger, J.-P., 1994.
\newblock A perfectly matched layer for absorption of electromagnetic waves,
  {\it Journal of Computational Physics\/}, {\bf 114}, 185--200.

\bibitem[Billette \& Brandsberg-Dahl(2004)]{Billette_2004_BPB}
Billette, F.~J. \& Brandsberg-Dahl, S., 2004.
\newblock The 2004 {BP} velocity benchmark, in {\em Extended Abstracts,
  67$^{th}$ Annual {EAGE} Conference \& Exhibition, Madrid, Spain\/}, p. B035.

\bibitem[B{\"o}hmer \& Stetter(1984)]{Bohmer_1984_DCM}
B{\"o}hmer, K. \& Stetter, H., 1984.
\newblock Defect correction methods, {\it Theory and applications\/}, {\bf 5}.

\bibitem[Boyd et~al.(2010)Boyd, Parikh, Chu, Peleato, \&
  Eckstein]{Boyd_2011_DOS}
Boyd, S., Parikh, N., Chu, E., Peleato, B., \& Eckstein, J., 2010.
\newblock Distributed optimization and statistical learning via the alternating
  direction of multipliers, {\it Foundations and trends in machine learning\/},
  {\bf 3}(1), 1--122.

\bibitem[Burvall et~al.(2011)Burvall, Lundstrom, Takman, Larsson, \&
  Hertz]{Burvall_2011_PRX}
Burvall, A., Lundstrom, U., Takman, P. A.~C., Larsson, D.~H., \& Hertz, H.~M.,
  2011.
\newblock Phase retrieval in {X}-ray phase-contrast imaging suitable for
  tomography, {\it Optics Express\/}, {\bf 19}(11), 10359--10376.

\bibitem[Candes et~al.(2013)Candes, Strohmer, \& Voroninski]{Candes_2013_PEA}
Candes, E.~J., Strohmer, T., \& Voroninski, V., 2013.
\newblock Phaselift: Exact and stable signal recovery from magnitude
  measurements via convex programming, {\it Communications on Pure and Applied
  Mathematics\/}, {\bf 66}(8), 1241--1274.

\bibitem[Chen et~al.(2013)Chen, Cheng, Feng, \& Wu]{Chen_2013_OFD}
Chen, Z., Cheng, D., Feng, W., \& Wu, T., 2013.
\newblock An optimal 9-point finite difference scheme for the {H}elmholtz
  equation with {PML}, {\it International Journal of Numerical Analysis \&
  Modeling\/}, {\bf 10}(2).

\bibitem[Choi \& Alkhalifah(2015)]{Choi_2015_UPI}
Choi, Y. \& Alkhalifah, T., 2015.
\newblock Unwrapped phase inversion with an exponential damping, {\it
  Geophysics\/}, {\bf 80}(5), 251--264.

\bibitem[Eldar et~al.(2016)Eldar, Hammen, \& Mixon]{Eldar_2016_RAP}
Eldar, Y.~C., Hammen, N., \& Mixon, D.~G., 2016.
\newblock Recent advances in phase retrieval, {\it IEEE signal processing
  magazine, Lecture Notes\/}, {\bf September}, 158--162.

\bibitem[Fienup(1982)]{Fienup_1982_PRA}
Fienup, J.~R., 1982.
\newblock Phase retrieval algorithms: a comparison, {\it Applied optics\/},
  {\bf 21}(15), 2758--2769.

\bibitem[Fogel et~al.(2016)Fogel, Waldspurger, \& d'Aspremont]{Fogel_2016_PRI}
Fogel, F., Waldspurger, I., \& d'Aspremont, A., 2016.
\newblock Phase retrieval for imaging problems, {\it Mathematical programming
  computation\/}, {\bf 8}(3), 311--335.

\bibitem[Fu \& Symes(2017)]{Fu_2017_DPM}
Fu, L. \& Symes, W.~W., 2017.
\newblock A discrepancy-based penalty method for extended waveform inversion,
  {\it Geophysics\/}, {\bf R282-R298}, 78--82.

\bibitem[Gerchberg(1972)]{Gerchberg_1972_APA}
Gerchberg, R.~W., 1972.
\newblock A practical algorithm for the determination of phase from image and
  diffraction plane pictures, {\it Optik\/}, {\bf 35}, 237--246.

\bibitem[Gholami(2014)]{Gholami_2014_PRT}
Gholami, A., 2014.
\newblock Phase retrieval through regularization for seismic problems, {\it
  Geophysics\/}, {\bf 79}(5), V153--V164.

\bibitem[Gholami \& Hosseini(2013)]{Gholami_2013_BCT}
Gholami, A. \& Hosseini, S.~M., 2013.
\newblock A balanced combination of {T}ikhonov and total variation
  regularizations for reconstruction of piecewise-smooth signals, {\it Signal
  Processing\/}, {\bf 93}, 1945--1960.

\bibitem[Gholami \& Naeini(2019)]{Gholami_2019_3DD}
Gholami, A. \& Naeini, E.~Z., 2019.
\newblock {3D} {D}ix inversion using bound-constrained {TV} regularization,
  {\it Geophysics\/}, {\bf 84}(3), 1--43.

\bibitem[Gholami et~al.(2018)Gholami, Aghamiry, \& Abbasi]{Gholami_2017_CNA}
Gholami, A., Aghamiry, H., \& Abbasi, M., 2018.
\newblock Constrained nonlinear {AVO} inversion using {Z}oeppritz equations,
  {\it Geophysics\/}, {\bf 83(3)}, R245--R255.

\bibitem[Goldstein \& Osher(2009)]{Goldstein_2009_SBM}
Goldstein, T. \& Osher, S., 2009.
\newblock The split {B}regman method for {L}1-regularized problems, {\it {SIAM}
  Journal on Imaging Sciences\/}, {\bf 2}(2), 323--343.

\bibitem[G{\'{o}}rszczyk et~al.(2017)G{\'{o}}rszczyk, Operto, \&
  Malinowski]{Gorszczyk_2017_TRW}
G{\'{o}}rszczyk, A., Operto, S., \& Malinowski, M., 2017.
\newblock Toward a robust workflow for deep crustal imaging by {FWI} of {OBS}
  data: The eastern nankai trough revisited, {\it Journal of Geophysical
  Research: Solid Earth\/}, {\bf 122}(6), 4601--4630.

\bibitem[Harrison(1993)]{Harrison_1993_PPI}
Harrison, R.~W., 1993.
\newblock Phase problem in crystallography, {\it JOSA a\/}, {\bf 10}(5),
  1046--1055.

\bibitem[Jiang et~al.(2017)Jiang, So, \& Liu]{Jiang_2017_RPR}
Jiang, X., So, H.~C., \& Liu, X., 2017.
\newblock Robust phase retrieval via {ADMM} with outliers, {\it
  arXiv:1702.06157v1\/}.

\bibitem[Kreutz-Delgado(2009)]{KREUTZ_2009_TCG}
Kreutz-Delgado, K., 2009.
\newblock The complex gradient operator and the cr-calculus, {\it arXiv
  preprint arXiv:0906.4835\/}.

\bibitem[Lange(2016)]{Lange_2016_MOA}
Lange, K., 2016.
\newblock {\it MM optimization algorithms\/}, vol. 147, SIAM.

\bibitem[Martin et~al.(2006)Martin, Wiley, \& Marfurt]{Martin_2006_M2E}
Martin, G.~S., Wiley, R., \& Marfurt, K.~J., 2006.
\newblock Marmousi2: An elastic upgrade for {M}armousi, {\it The Leading
  Edge\/}, {\bf 25}(2), 156--166.

\bibitem[Meyer(2000)]{Meyer_2000_MAA}
Meyer, C.~D., 2000.
\newblock {\it Matrix analysis and applied linear algebra\/}, vol.~71, Siam.

\bibitem[Millane(1990)]{Millane_1990_PRI}
Millane, R.~P., 1990.
\newblock Phase retrieval in crystallography and optics, {\it JOSA A\/}, {\bf
  7}(3), 394--411.

\bibitem[Netrapalli et~al.(2013)Netrapalli, Jain, \&
  Sanghavi]{Netrapalli_2013_PRU}
Netrapalli, P., Jain, P., \& Sanghavi, S., 2013.
\newblock Phase retrieval using alternating minimization, in {\em Advances in
  Neural Information Processing Systems\/}, pp. 2796--2804.

\bibitem[Nocedal \& Wright(2006)]{Nocedal_2006_NO}
Nocedal, J. \& Wright, S.~J., 2006.
\newblock {\it Numerical Optimization\/}, Springer, 2nd edn.

\bibitem[Oppenheim \& Lim(1981)]{Oppenheim_1981_IPS}
Oppenheim, A.~V. \& Lim, J.~S., 1981.
\newblock The importance of phase in signals, {\it Proceedings of the IEEE\/},
  {\bf 69}(5), 529--541.

\bibitem[Pfeiffer et~al.(2006)Pfeiffer, Weitkamp, Bunk, \&
  David]{Pfeiffer_2006_PRD}
Pfeiffer, F., Weitkamp, T., Bunk, O., \& David, C., 2006.
\newblock Phase retrieval and differential phase-contrast imaging with
  low-brillance {X}-ray sources, {\it Nature Physics\/}, {\bf 2}, 258--261.

\bibitem[Plessix(2006)]{Plessix_2006_RAS}
Plessix, R.~E., 2006.
\newblock A review of the adjoint-state method for computing the gradient of a
  functional with geophysical applications, {\it Geophysical Journal
  International\/}, {\bf 167}(2), 495--503.

\bibitem[Pratt et~al.(1998)Pratt, Shin, \& Hicks]{Pratt_1998_GNF}
Pratt, R.~G., Shin, C., \& Hicks, G.~J., 1998.
\newblock {G}auss-{N}ewton and full {N}ewton methods in frequency-space seismic
  waveform inversion, {\it Geophysical Journal International\/}, {\bf 133},
  341--362.

\bibitem[Qian et~al.(2017)Qian, Fu, Sidiropoulos, Huang, \& Xie]{Qian_2017_IAO}
Qian, C., Fu, X., Sidiropoulos, N.~D., Huang, L., \& Xie, J., 2017.
\newblock Inexact alternating optimization for phase retrieval in the presence
  of outliers, {\it IEEE Transactions on Signal Processing\/}, {\bf 65}(22),
  6069--6082.

\bibitem[Shah(2014)]{Shah_PhD_2014}
Shah, N.~K., 2014.
\newblock {\it Seismic Full Waveform Inversion for wrapped and unwrapped
  phase\/}, Ph.D. thesis, Imperial College London.

\bibitem[Shechtman et~al.(2015)Shechtman, Eldar, Cohen, Chapman, Miao, \&
  Segev]{Shechtman_2015_PRW}
Shechtman, Y., Eldar, Y.~C., Cohen, O., Chapman, H.~N., Miao, J., \& Segev, M.,
  2015.
\newblock Phase retrieval with application to optical imaging: a contemporary
  overview, {\it IEEE signal processing magazine\/}, {\bf 32}(3), 87--109.

\bibitem[Shipp \& Singh(2002)]{Shipp_2002_TDF}
Shipp, R.~M. \& Singh, S.~C., 2002.
\newblock Two-dimensional full wavefield inversion of wide-aperture marine
  seismic streamer data, {\it Geophysical Journal International\/}, {\bf 151},
  325--344.

\bibitem[Sirgue \& Pratt(2004)]{Sirgue_2004_EWI}
Sirgue, L. \& Pratt, R.~G., 2004.
\newblock Efficient waveform inversion and imaging : a strategy for selecting
  temporal frequencies, {\it Geophysics\/}, {\bf 69}(1), 231--248.

\bibitem[{van Leeuwen} \& Herrmann(2016)]{vanLeeuwen_2016_PMP}
{van Leeuwen}, T. \& Herrmann, F., 2016.
\newblock A penalty method for {PDE}-constrained optimization in inverse
  problems, {\it Inverse Problems\/}, {\bf 32(1)}, 1--26.

\bibitem[{van Leeuwen} \& Herrmann(2013)]{VanLeeuwen_2013_MLM}
{van Leeuwen}, T. \& Herrmann, F.~J., 2013.
\newblock Mitigating local minima in full-waveform inversion by expanding the
  search space, {\it Geophysical Journal International\/}, {\bf 195(1)},
  661--667.

\bibitem[Virieux \& Operto(2009)]{Virieux_2009_OFW}
Virieux, J. \& Operto, S., 2009.
\newblock An overview of full waveform inversion in exploration geophysics,
  {\it Geophysics\/}, {\bf 74}(6), WCC1--WCC26.

\bibitem[Waldspurger et~al.(2015)Waldspurger, d’Aspremont, \&
  Mallat]{Waldspurger_2015_PRM}
Waldspurger, I., d’Aspremont, A., \& Mallat, S., 2015.
\newblock Phase recovery, maxcut and complex semidefinite programming, {\it
  Mathematical Programming\/}, {\bf 149}(1-2), 47--81.

\bibitem[Walther(1963)]{Walther_1963_TQO}
Walther, A., 1963.
\newblock The question of phase retrieval in optics, {\it Optica Acta:
  International Journal of Optics\/}, {\bf 10}(1), 41--49.

\bibitem[Wu et~al.(1983)]{Wu_1983_OTC}
Wu, C.~J. et~al., 1983.
\newblock On the convergence properties of the em algorithm, {\it The Annals of
  statistics\/}, {\bf 11}(1), 95--103.

\end{thebibliography}

\newcommand{\SortNoop}[1]{}

\begin{figure}
\center
      \includegraphics[trim={0.1cm 0.7cm 0 0},clip,width=0.48\textwidth]{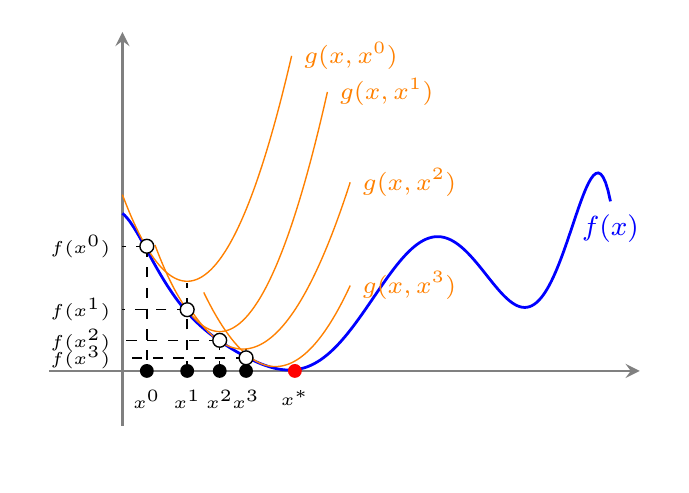} 
\caption{Sketch of the iterative MM algorithm to solve optimization problem \ref{c0_a}. The function $f(\bold{x})$ is shown in blue while a few surrogate functions ${g(\bold{x},\bold{x}^k)}, k \in \{0, 1, 2, 3\}$ are shown in orange. The MM algorithm seeks to find a local minimizer of $f(\bold{x})$, red point, by iteratively minimizing easy to minimize functions ${g(\bold{x},\bold{x}^k)}$. The surrogate function ${g(\bold{x},\bold{x}^k)}$ is greater than or equal to $f(\bold{x})$, and the equality holds at the current optimal point $\bold{x}^k$. Minimization of ${g(\bold{x},\bold{x}^k)}$ gives a new optimal point and it reaches to $\bold{x}^*$ when $k$ tends to infinity.}
\label{fig:su}
\end{figure}

%
\append{Majorization-minimization to solve phase retrieval problem} \label{Appa}
Here, we present a simple algorithm to solve phase retrieval problem based on the majorization-minimization (MM) \citep{Lange_2016_MOA} technique.
The governing idea of MM is to find the minimum of non-convex/convex function $f(\bold{x})$ via the iterative minimization of a simpler convex surrogate function $g(\bold{x},\bold{x}^k)$ that majorizes $f(\bold{x})$ at step $k$ (i.e. $g(\bold{x},\bold{x}^k) \geq f(\bold{x})$).
Figure \ref{fig:su} shows a schematic of the MM process. The non-convex function ${f(\bold{x})}$ is shown in blue while a few surrogate functions ${g(\bold{x},\bold{x}^k)}$ for points ${\bold{x}}^k$, $k\in \{0, 1, 2 ,3\}$, are shown in orange. This figure shows how the iterative MM algorithm approaches a local minimum of $f(\bold{x})$ through the minimum of easy to minimize surrogate functions ${g(\bold{x},\bold{x}^k)}$.\\   
The phase-retrieval problem for linear system $\bold{Lx}=\bold{y}$ can be written as \citep{Gholami_2014_PRT} 
 \begin{equation} \label{c0_a}
 \min_{\bold{x}} ~~ f(\bold{x})=\min_{\bold{x}}~~\frac12 \|\bold{|Lx|}-\bold{|y|}\|_2^2,
 \end{equation}
where $\bold{L}\in \mathbb{C}^{m\times n}$,  $\bold{x}\in \mathbb{C}^{n\times 1}$, and  $\bold{y} \in \mathbb{C}^{m\times 1}$.
Problem \eqref{c0_a} is non-convex where this non-convexity finds its root in removing the phase of the right hand side.  

Let us introduce the auxiliary variable $\bold{z=Lx}$, then the objective function in \eqref{c0_a} can be written as
\begin{equation} \label{c0_aa}
\tilde{f}(\bold{z})=\frac12 \|\sqrt{\bold{z}{\bold{z}^*}}-\bold{|y|}\|_2^2. 
\end{equation}

Using the Wirtinger calculus \citep{KREUTZ_2009_TCG}, 
the quadratic approximation of $f(\bold{x})$ around an initial guess $\bold{x}^k$ is
\begin{eqnarray} \label{QUAD}
g(\bold{x},\bold{x}^k) = f(\bold{x}^k) + \nabla f(\bold{x}^k)^T\Delta\bold{x}
+\Delta\bold{x}^T\bold{H}(\bold{x}^k)\Delta\bold{x},
\end{eqnarray}
where $\Delta\bold{x}=\bold{x-x}^k$, $\nabla f$ is the gradient of $f(\bold{x})$ and $\bold{H}$ denotes the Hessian, both evaluated at $\bold{x}^k$. 
%
 Simple algebra shows that the gradient of $f(\bold{x})$ is given by 
 \begin{align}  \label{grad}
 \nabla f(\bold{x}) &= 2 \bold{L}^T \frac{\partial \tilde{f}(\bold{z})}{\partial {\bold{z}^*}} \\
 &= \bold{L}^T \text{diag}({\bf 1} - \frac{|\bold{y}|}{|\bold{Lx}|})\bold{Lx} \\
&= \bold{L}^T (\bold{Lx} - |\bold{y}|e^{j\angle \bold{Lx}}), \label{grad1}
\end{align}
where $\bf{1}$ is an all-ones vector, and the Hessian operator is given by
\begin{eqnarray}
\label{hessian}
\bold{H}(\bold{x}) &=& 2\bold{L}^T \left( \frac{\partial^2 \tilde{f}(\bold{z}) }{\partial\bold{z}\partial\bold{z}^*}  
 +   \frac{\partial^2 \tilde{f}(\bold{z}) }{\partial\bold{z}^*\partial\bold{z}^*}  \boldsymbol{\Omega} \right) \bold{L}\\
&=& \bold{L}^T \left(\text{diag}({\bf 1} - \frac{|\bold{y}|}{2|\bold{Lx}|}) + \text{diag}(\frac{|\bold{y}|e^{j2\angle \bold{Lx}}}{2|\bold{Lx}|}) \boldsymbol{\Omega}\right) \bold{L}, \nonumber 
\end{eqnarray}
where $\boldsymbol{\Omega}(\bold{z})=\bold{z}^*$ denotes the complex conjugate operator.

Ignoring the nonlinear terms of Hessian leads to the following Gauss-Newton approximation
\begin{eqnarray} \label{QUASI}
g(\bold{x},\bold{x}^k) = f(\bold{x}^k) + \nabla f(\bold{x}^k)^T\Delta\bold{x}
+\Delta\bold{x}^T\bold{L}^T\bold{L}\Delta\bold{x},
\end{eqnarray}
or equivalently
%
\begin{align} \label{quad}
g(\bold{x},\bold{x}^k)=\frac12 \|\bold{Lx}-\bold{|y|}e^{j\angle \bold{Lx}^k}\|_2^2.
\end{align} 
The quadratic function $g(\bold{x},\bold{x}^k)$ majorizes $f(\bold{x})$ at the point $\bold{x}^k$ provided that the following MM conditions are satisfied: 
\begin{subequations}
    \begin{empheq}[left={\empheqlbrace\,}]{align}
      & g(\bold{x}^k,\bold{x}^k)  =  f(\bold{x}^k) \label{cond1} \\
      & g(\bold{x},\bold{x}^k)  \geq  f(\bold{x})~~~~~~~ \text{for all}~ \bold{x}. \label{cond2} 
    \end{empheq}
\end{subequations}
Equality \eqref{cond1} can easily be confirmed by substituting $\bold{x}$ by $\bold{x}^k$ in $f(\bold{x})$ and $g(\bold{x},\bold{x}^k)$.
Furthermore, let $z_i=(\bold{Lx})_i$ and $\theta_i=\angle(\bold{Lx}^k)_i$ then according to the backward triangle inequality \citep{Meyer_2000_MAA} 
\begin{align}
\frac12 |z_i-|y_i|e^{j\theta_i}|^2 &\geq \frac12 (|z_i|-||y_i|e^{j\theta_i}|)^2 \\
&= \frac12 (|z_i|-|y_i||e^{j\theta_i}|)^2 \nonumber \\
&= \frac12 (|z_i|-|y_i|)^2. \nonumber 
\end{align}
Hence,
\begin{equation}
\frac12\sum_i |z_i-|y_i|e^{j\theta_i}|^2\geq  \frac12 \sum_i (|z_i|-|y_i|)^2,
\end{equation}
which confirms that the inequality \eqref{cond2} is also satisfied and thus $g(\bold{x},\bold{x}^k)$ majorizes $f(\bold{x})$.  
These conditions guarantee that iterative minimization of $g(\bold{x},\bold{x}^k)$ converge to a local minimum of $f(\bold{x})$ as $k$  tends to infinity \citep{Wu_1983_OTC}.\\

\end{document}